\documentclass{siamart251216}
\usepackage{verbatim}
\usepackage{amsmath,amsfonts}
\usepackage{amssymb}
\usepackage{graphicx}
\usepackage{wasysym}
\newtheorem{thm}{Theorem}
\newtheorem{remark}{Remark}
\newtheorem{cor}{Corollary}
\newtheorem{prop}{Proposition}
\newtheorem{example}{Example}
\newtheorem{lem}{Lemma}
\newtheorem{defn}{Definition}

\begin{document}
\title{Dynamical Unbalanced Optimal Transport: Duality, Reachability, and
Computation}

\author{Dongjun Wu 
	\thanks{D. Wu is with the Department of Electrical Engineering, KU Leuven, Belgium}
	\and
	Anders Rantzer
	\thanks{A. Rantzer is with the Department of Automatic Control, Lund University, Sweden}
}

\maketitle

\begin{abstract}
We study dynamical unbalanced optimal transport for nonlinear control-affine
systems in which the control input affects both the transportation
of mass and its growth or decay. This phenomenon arises naturally
in controlled transport-reaction systems such as cancer treatment,
and distinguishes the problem from standard unbalanced optimal transport,
where transport and mass variation are typically independent. To address
this problem, we combine classical relaxed-control theory with a weighted
superposition principle to obtain an infinite-dimensional linear programming
reformulation over weighted trajectories. Under suitable convexity
and compactness assumptions, we show that this relaxation is tight.
Based on the relaxation, we derive a Kantorovich dual formulation
in terms of Hamilton-Jacobi-Bellman inequalities, and establish strong
duality over smooth subsolutions. We then utilize the dual formulation
to characterize reachability of the system, namely, when a given terminal
measure is reachable from a given initial measure. In addition, tractable
sufficient conditions for reachability are provided. Finally, we develop
an augmented-Lagrangian algorithm based on the dual formulation for
numerical computation of the optimal transport plan. The framework
is illustrated on a phenotype-structured cancer treatment problem
in which therapy simultaneously drives phenotypic adaptation and population
growth or death.
\end{abstract}

\section{Introduction}

\label{sec:intro}

Optimal transport (OT) provides a variational framework for transforming
one distribution into another while minimizing a transportation cost.
Its static and dynamic formulations have become important tools in
applied mathematics
\cite{villaniOptimalTransport2009,santambrogioOptimalTransport2015,peyreCuturiComputational2019}.
In particular, the dynamical formulation introduced by Benamou and Brenier \cite{benamou2000computational} 
has found rich applications in systems and control
\cite{chenOptimalTransportSystems2021,peyreCuturiComputational2019}. Indeed, the Benamou-Brenier problem can
be seen as minimizing the cost of transporting a distribution of controlled particles.

In the original Benamou-Brenier problem, each particle follows a simple first-order integrator dynamics.
Recent years have seen extensions to more general dynamics with possibly state and input constraints 
\cite{chenOptimalTransportSystems2021,bernardOptimalMassTransportation2007,fornasierMeanFieldGamma2019,cavagnariLagrangianEulerianKantorovich2022,ghoussoubOptimalTransportControlled2018,elamvazhuthiDynamicalOptimalTransport2023}.
The subject has become what is now called \emph{dynamical optimal transport}.
In particular, optimal transport over linear and affine
nonlinear dynamics has been extensively studied \cite{Chen2017linear,elamvazhuthiDynamicalOptimalTransport2023}.
These works are, however, intrinsically mass preserving, which limits
its applicability to systems in which agents, material, probability
mass, or other resources may be created or removed during the evolution. 

Several extensions of OT have therefore been developed with non-conservation
of masses. These include generalized Wasserstein metrics for transport
equations with sources \cite{piccoliRossiGeneralizedWasserstein2014,kondratyevNewOptimalTransport2016},
the Wasserstein-Fisher-Rao/Hellinger-Kantorovich metrics \cite{chizatInterpolatingDistanceOptimal2018,lieroOptimalTransportReaction2016},
and dynamic and entropy-transport formulations of unbalanced optimal
transport (UOT) \cite{chizatUnbalancedOptimalTransport2018,liero2018optimal}.
Relaxation, duality, optimality, and Monge-type formulations have
subsequently been developed in considerable generality \cite{savareSodiniRelaxation2024}.
Partial transport provides another mechanism for handling unequal
masses, where only a prescribed fraction of mass is transported while
the remaining mass is discarded \cite{caffarelliMcCannPartialTransport2010,figalliPartialTransport2010}. 

An important feature of UOT problem is that the transport velocity
and the rate of mass-variation are treated as independent optimization
variables \cite{chizatUnbalancedOptimalTransport2018}. In controlled
systems, however, these two mechanisms may be coupled through the
same actuator. For example, a control input may simultaneously affect
the motion of an individual particle and its rate of growth or decay.
This coupling is not merely a mathematical generalization. For instance,
in phenotype-structured models of cancer treatment, the drug intensity
may affect both the evolution of cellular phenotypes and the net proliferation
or death rate of the population \cite{chisholmDrugTolerance2015,lorzPhenotypicHeterogeneity2015,almeidaCytotoxicTherapy2019}.
At the level of an individual cell or subpopulation, treatment can
drive phenotypic adaptation toward drug-tolerant states while at the
same time suppressing proliferation or inducing cell death. Thus the
transport of mass in phenotype space and the creation or destruction
of mass are controlled by the same input. 

To the best of our knowledge, the corresponding problem in which transport
velocity and local mass variation are coupled has not been developed,
apart from the preliminary conference version of the present work
\cite{wuRantzerCDC2025}. 

To approach the problem, we employ the following controlled transport-reaction
model in which the reaction term on the right accounts for mass variation:
\begin{equation}
\partial_{t}\mu+\nabla\cdot\bigl(F(t,x,u)\mu\bigr)=G(t,x,u)\mu.\label{eq:intro-dynamics}
\end{equation}
Both the transport field $F$ and the growth rate $G$ are nonlinear
in the state and affine in the control variable $u$. To see how the
total mass evolves, suppose $\mu_t\in\mathcal{M}_{+}(\mathbb{R}^{n})$
and let $m(t)=\int d\mu_{t}$ be the total mass at time $t$. Then,
formally, 
\[
\dot{m}(t)=\int G(t,x,u)d\mu_{t}.
\]
Hence $G$ governs local creation and destruction of mass. The salient
feature of \eqref{eq:intro-dynamics} is that the velocity $F$ and
the reaction rate $G$ cannot, in general, be selected independently,
since both are influenced by the same control input $u$, manifesting
a key difference compared to standard UOT.

Given the two boundary measures $\mu_{0}$ and $\mu_{T}$, our objective
is to steer $\mu_{0}$ to $\mu_{T}$ while minimizing a Lagrangian type cost functional,
\[
\int^{T}_{0}\int L(t,x,u(t,x))d\mu_{t}(x)dt.
\]
The resulting optimization problem is nonlinear. To expose its
convex structure, we combine classical relaxed-control theory
\cite{young1969lectures} with a
weighted superposition principle \cite{brediesSuperpositionInhomogeneous2022} for the equation
\eqref{eq:intro-dynamics}.
This leads to a path-space
formulation over \emph{weighted relaxed trajectories}, and hence to an
infinite-dimensional linear program. Under suitable convexity and compactness  
assumptions imposed, the relaxation is tight, namely, the
path-space problem has the same value as the original optimization
problem. 

The relaxation result provides the foundation for the dual formulation developed subsequently.
Once the problem is lifted to the weighted path space, convex duality reduces the endpoint-constrained transport
problem to a family of single-trajectory optimal control problems. Their value functions satisfy a
Hamilton-Jacobi-Bellman (HJB) inequality.
These value functions, however, need not be smooth, so convex duality alone does not justify restricting the
dual problem to classical HJB subsolutions.
In optimal control, this issue was solved by Fleming-Vermes and Vinter, who established dual formulations
over classical $C^{1}$ subsolutions for nonlinear optimal control problems
\cite{vinterConvexDualityNonlinear1993,flemingConvexDualityApproach1989}.
In parallel, HJB duality over smooth subsolutions have also been developed in balanced OT 
\cite{bernardOptimalMassTransportation2007,ghoussoubOptimalTransportControlled2018}, and UOT literature
\cite{chizatThesis2017}. 
However, for
the general dynamics \eqref{eq:intro-dynamics}, analogous dual formulation over classical smooth HJB subsolutions
has not been developed; the main difficulty lies in the lack of coercivity. 
Following
\cite{vinterConvexDualityNonlinear1993}, we establish such a result. This yields strong duality between the
original problem, and the smooth HJB dual.

An important question regarding UOT of system \eqref{eq:intro-dynamics} is reachability.
That is, given $\mu_0$ and an admissible control input set, when is a measure $\mu_T$ reachable?
Our duality result provides a mechanism to address this question. 
In the literature, reachability in spaces of probability measures have been studied
through controlled continuity equations, differential inclusions,
and HJB methods \cite{cavagnariGeneralizedControlSystems2018,cavagnariAveragedTimeOptimal2018,cavagnariMarigondaAttainability2021}.
Recent results also relate reachability of a mass-preserving continuity
equation to controllability of an underlying control-affine system
\cite{elamvazhuthiDynamicalOptimalTransport2023},
and to realizability of transport maps \cite{caponigroTransportMaps2025}.
The present setting differs from these works in that the state of
interest is a positive measure whose total mass may vary, while its
transport and mass evolution are coupled through the same control.
We characterize the resulting reachable set through the duality theory
that we developed and propose tractable sufficient conditions for
reachability.

There is also a substantial computational literature for dynamic OT
and UOT. The Benamou-Brenier discretization has been combined with
proximal splitting and first-order convex optimization
\cite{papadakisOptimalTransportProximal2014,peyreCuturiComputational2019}.
For UOT, entropic scaling algorithms and unbalanced Sinkhorn methods
provide fast alternatives in static and regularized settings
\cite{chizatScalingAlgorithms2018,sejourneSinkhornUnbalanced2019}, while
Wasserstein-Fisher-Rao splitting schemes exploit the decomposition into
transport and reaction steps \cite{gallouetJKO2017,gallouetLabordeMonsaingeon2019}.
Algorithm for the general dynamical UOT under dynamics \eqref{eq:intro-dynamics} 
has not been developed.

A preliminary conference version \cite{wuRantzerCDC2025} of the present work introduced the
controlled dynamical UOT model, a Vinter type relaxation, a preliminary duality result (with some key proof
steps left out), and a
first-order augmented-Lagrangian scheme. The present paper substantially
extends that work. First, the relaxed formulation is developed as a
weighted path-space problem derived from a superposition principle and
is proved to be exactly equivalent to the Eulerian problem. Second, the
restriction of the dual to $C^{1}$ functions in the dual formulation is now proven rigorously by a 
smooth-subsolution approximation theorem. Third, we develop a 
reachability theory based on duality, going beyond the sufficient controllability
condition in the conference version. Finally, Hamiltonian conjugacy is
used to identify the augmented-Lagrangian multipliers with the discrete
mass, transport flux, and reaction flux, and the numerical study is
expanded to a nonlinear transport-reaction model. The main
\textbf{contributions} are as follows.
\begin{enumerate}
\item We formulate dynamical unbalanced optimal transport for nonlinear
control-affine systems in which the same control input affects both
particle motion and mass variation. By combining relaxed-control theory
with the weighted superposition principle, we obtain an
infinite-dimensional linear program over weighted relaxed trajectories,
prove exactness of the relaxation, and establish strong duality with an
HJB formulation over $C^{1}$ subsolutions.
\item We develop a reachability theory for the dynamical unbalanced optimal
transport problem. The reachable set is characterized through its
support function and HJB inequalities, yielding smooth separating
certificates for unreachable target measures. We also relate measure
reachability to controllability of an augmented finite-dimensional
system and obtain a direct characterization in the balanced case.
\item 
We develop a splitting algorithm for the dynamical UOT problem based on the 
dual formulation in which the nonlinear constraint is enforced by
pointwise projections and the dual-potential update reduces to a sparse
linear system. A new primal-dual interpretation shows directly from
Hamiltonian conjugacy that the multiplier recovers the discrete mass,
transport flux, and reaction flux, and hence the control input. The method is illustrated on a
nonlinear phenotype-structured cancer-treatment model.
\end{enumerate}
The remainder of the paper is organized as follows. Section~\ref{sec:problem}
introduces the controlled dynamical unbalanced optimal transport problem.
Section~\ref{sec:relaxation} develops the weighted relaxed-arc formulation
and proves its equivalence with the original problem. Section~\ref{sec:duality}
establishes the HJB duality and the smooth subsolution approximation.
Section~\ref{sec:reachability} studies reachability of positive measures.
Section~\ref{sec:numerics} develops the augmented-Lagrangian discretization
and the recovery of primal variables and feedback.

\emph{Notations.} $C_{b}(X)$ is the set of real-valued bounded continuous
functions on $X$. $C^{r}_{c}(X)$ is the space of $C^{r}$ functions
with compact support in $X$. $\mathcal{M}_{+}(X)$ is the set of
non-negative Borel measures on $X$. Let $F:X\to Y$ be a measurable mapping and $\mu$ be a measure on $X$; 
the pushforward measure
$F_\# \mu$ is a measure on $Y$ defined as $(F_\# \mu)(A) = \mu(F^{-1}(A))$. Let $r>0$, then
$\overline{B}_r$ is the closed ball with radius $r$.

\section{Problem formulation } \label{sec:problem}

Consider the controlled transport-reaction equation 
\begin{equation}
\partial_{t}\mu+\nabla\cdot\bigl(F(t,x,u)\mu\bigr)=G(t,x,u)\mu,\label{sys:PDE}
\end{equation}
where $t\in[0,T]$, $x\in\mathbb{R}^{n}$, $u\in U\subset\mathbb{R}^{m}$,
and 
\begin{equation}
\begin{aligned}F(t,x,u) & =f_{0}(t,x)+\sum^{m}_{i=1}u_{i}f_{i}(t,x),\\
G(t,x,u) & =g_{0}(t,x)+\sum^{m}_{i=1}u_{i}g_{i}(t,x).
\end{aligned}
\label{eq:F-G}
\end{equation}
Thus $F$ and $G$ are affine in the control variable and may be nonlinear
in the state variable. A measure curve $\mu:[0,T]\to\mathcal{M}_{+}(\mathbb{R}^{n})$
is said to solve \eqref{sys:PDE} in the distributional sense if,
for every $\varphi\in C^{1}_{c}([0,T]\times\mathbb{R}^{n})$, 
\begin{equation}
\begin{aligned} & \int^{T}_{0}\int_{\mathbb{R}^{n}}\Bigl(\partial_{t}\varphi+\nabla\varphi^{\top}F(t,x,u)+G(t,x,u)\varphi\Bigr)d\mu_{t}dt\\
 & \qquad\qquad\qquad\quad=\int_{\mathbb{R}^{n}}\varphi(T,\cdot)d\mu_{T}-\int_{\mathbb{R}^{n}}\varphi(0,\cdot)d\mu_{0}.
\end{aligned}
\label{eq:weak-pde}
\end{equation}
We say that $\mu$ is \emph{narrowly continuous} if, for every $\psi\in C_{b}(\mathbb{R}^{n})$,
the mapping 
\[
t\longmapsto\int_{\mathbb{R}^{n}}\psi(x)\,d\mu_{t}(x)
\]
is continuous. The space of such curves is denoted by $C_{w}([0,T];\mathcal{M}_{+}(\mathbb{R}^{n}))$.

\medskip{}
 \noindent\textbf{Assumption A.} The following hold. 
\begin{enumerate}
\item $U\subset\mathbb{R}^{m}$ is nonempty, compact, and convex. 
\item The vector fields $f_{i}:[0,T]\times\mathbb{R}^{n}\to\mathbb{R}^{n}$
and functions $g_{i}:[0,T]\times\mathbb{R}^{n}\to\mathbb{R}$ are
continuous and locally Lipschitz in $x$, uniformly in $t$, and $F$
satisfies the uniform linear-growth condition 
\begin{equation}
\sup_{u\in U}|F(t,x,u)|\le C(1+|x|),\label{eq:linear-growth}
\end{equation}
for some $C>0$ and all $(t,x)\in[0,T]\times\mathbb{R}^{n}$. 
\item The running cost $L:[0,T]\times\mathbb{R}^{n}\times U\to[0,+\infty)$
is continuous, locally Lipschitz in $x$ uniformly in $(t,u)$ on
compact sets, and $u\mapsto L(t,x,u)$ is convex for every $(t,x)$. 
\end{enumerate}
We call $u:[0,T]\times\mathbb{R}^{n}\to U$ an \emph{admissible control}
if it is measurable in $t$ and locally Lipschitz in $x$, uniformly
in $t$ on compact subsets of $\mathbb{R}^{n}$. The well-posedness
of the equation \eqref{sys:PDE} is guaranteed by the proposition
below. 
\begin{prop}[Well-posedness]
\label{prop:sol-PDE} Suppose Assumption A holds, let $u$ be an
admissible control, and let $\mu_{0}\in\mathcal{M}_{+}(\mathbb{R}^{n})$
have compact support. Let $\phi_{s,t}$ denote the flow generated
by 
\begin{equation}
\dot{x}(t)=F(t,x(t),u(t,x(t))).\label{eq:closed-loop-flow}
\end{equation}
Then \eqref{sys:PDE} admits a unique narrowly continuous distributional
solution, given by 
\begin{equation}
\mu_{t}=(\phi_{0,t})_{\#}\bigl(h(t,\cdot)\mu_{0}\bigr),\label{eq:phi-mu}
\end{equation}
where 
\begin{equation}
h(t,x)=\exp\left(\int^{t}_{0}G\bigl(s,\phi_{0,s}(x),u(s,\phi_{0,s}(x))\bigr)\,ds\right).\label{eq:weight-h}
\end{equation}
Moreover, if $K\subset\mathbb{R}^{n}$ is compact and ${\rm supp}\mu_{0}\subset K$,
then 
\begin{equation}
{\rm supp}\mu_{t}\subset\overline{B}_{R},\quad\forall t\in[0,T],\label{eq:compact-support-classical}
\end{equation}
where $R$ is a constant dependent on $R$ and $K$. 
\end{prop}
\begin{proof}
The representation and uniqueness are standard for linear non-homoge\-neous
transport equations with locally Lipschitz velocity fields; see \cite{ambrosioGigliSavare2008,maniglia2007}.
The linear-growth bound gives 
\[
1+|\phi_{s,t}(x)|\le(1+|x|)e^{C(t-s)}
\]
by Gronwall's inequality. Since $h(t,x)>0$, the representation \eqref{eq:phi-mu}
yields $\operatorname{supp}\mu_{t}\subset\phi_{0,t}(\operatorname{supp}\mu_{0})$,
and \eqref{eq:compact-support-classical} follows. 
\end{proof}

By Proposition \ref{prop:sol-PDE}, every solution of \eqref{sys:PDE}
issued from $K$ is supported in the common compact set $\overline{B}_{R}$.
Throughout the paper, all initial measures are supported in $K$ and
hence all measure curves may be regarded as curves in $\mathcal{M}_{+}(X)$,
with $X=\overline{B}_{R}$. If ${\rm supp}\mu_{T}\not\subset X$,
then the endpoint problem is infeasible. From now on, we assume 
\[
{\rm supp}\mu_{0}\subset K,\quad{\rm supp}\mu_{T}\subset X
\]
In particular, $F$, $G$, and $L$ are bounded on $[0,T]\times X\times U$
whenever needed.
\begin{defn}[Admissible Eulerian pair]
\label{def:eulerian-admissible} A pair $(\mu,u)$ is called an\emph{
admissible Eulerian pair} if $u:[0,T]\times\mathbb{R}^{n}\to U$ is
Borel measurable, $\mu\in C_{w}([0,T];\mathcal{M}_{+}(\mathbb{R}^{n}))$
satisfies \eqref{sys:PDE} in the distributional sense, and 
\begin{equation}
\int^{T}_{0}\int_{\mathbb{R}^{n}}\bigl(|F(t,x,u)|^{2}+|G(t,x,u)|^{2}\bigr)d\mu_{t}dt<\infty.\label{eq:HK-energy}
\end{equation}
\end{defn}

The additional integrability condition is precisely the one needed
for the weighted superposition theorem used below. Every regular control
with compactly supported initial measure generates an admissible pair
by Proposition \ref{prop:sol-PDE}.

In Definition \ref{def:eulerian-admissible}, we have enlarged the
admissible control input actions , because the convex relaxation developed
in the next section naturally produces Borel feedbacks rather than
locally Lipschitz ones. 

Now we are ready to define the \emph{Eulerian primal problem}. Fix
compactly supported endpoint measures $\mu_{0}\in\mathcal{M}_{+}(K)$,
$\mu_{T}\in\mathcal{M}_{+}(X)$. For an admissible pair, define 
\begin{equation}
J_{E}(\mu,u):=\int^{T}_{0}\int_{\mathbb{R}^{n}}L(t,x,u(t,x))d\mu_{t}dt.\label{eq:Euler-cost}
\end{equation}

\medskip{}
 \noindent\textbf{Problem E.} Let $\mu_{0}\in\mathcal{M}_{+}(K),\,\mu_{T}\in\mathcal{M}_{+}(X)$,
define
\begin{equation}
V_{E}(\mu_{0},\mu_{T})=\inf\bigl\{ J_{E}(\mu,u):\mu(0)=\mu_{0},\,\mu(T)=\mu_{T}\bigr\}\label{eq:VE}
\end{equation}
where infimum ranges over admissible pair $(\mu,u)$. If no admissible
pair connects $\mu_{0}$ to $\mu_{T}$, set $V_{E}(\mu_{0},\mu_{T})=+\infty$.

\section{Relaxation by weighted relaxed arcs}

\label{sec:relaxation}

We now derive a Lagrangian convex representation of Problem E in path
space. This formulation is particularly useful for endpoint duality and
reachability. The construction combines the classical relaxed-control
idea of Young and Vinter with the weighted superposition principle for
transport-reaction equations. 

\subsection{Weighted superposition}

Consider the following space of curves
\[
\mathcal{A}_{0}=AC([0,T];X)\times AC([0,T];\mathbb{R}_{>0}).
\]
Here $AC([0,T];X)$ represents the set of absolutely continuous curves
in $X$ on the interval $[0,T]$. An element in $\mathcal{A}_{0}$
is denoted $(\gamma,h)$. For any curve $t\mapsto c(t)$, let $e_{t}$
denote the evaluation map at $t$, i.e., $e_{t}c=c(t)$. For $z=(\gamma,h)\in\mathcal{A}_{0}$,
let $\pi_{1}$ and $\pi_{2}$ be the projections over the first and
second entry respectively, i.e, $\pi_{1}z=\gamma$ and $\pi_{2}z=h$. 

Our relaxation relies on the following weighted superposition principle
essentially proven in \cite{brediesSuperpositionInhomogeneous2022}. 
\begin{prop}[Weighted superposition]
\label{prop:BCF-controlled} Let $(\mu,u)$ be an admissible Eulerian
pair with ${\rm supp}\mu_{0}\subset K$. Then there exists a finite
nonnegative measure $\sigma\in\mathcal{M}_{+}(\mathcal{A}_{0})$ such
that 
\begin{equation}
(e_{0}\circ\pi_{1})_{\#}\sigma=\mu_{0},\quad h(0)=1\quad\sigma\text{-a.e.},\label{eq:path-initial}
\end{equation}
and, for every $t\in[0,T]$, 
\begin{equation}
\mu_{t}=\int h(t)\delta_{\gamma(t)}d\sigma(\gamma,h).\label{eq:weighted-superposition}
\end{equation}
Moreover, $\sigma$ is concentrated on curves satisfying 
\begin{align}
\dot{\gamma}(t) & =F(t,\gamma(t),u(t,\gamma(t))),\label{eq:path-x}\\
\dot{h}(t) & =G(t,\gamma(t),u(t,\gamma(t)))h(t)\label{eq:path-h}
\end{align}
for a.e. $t\in[0,T]$. 
\end{prop}
\begin{proof}
Apply Theorem 1.1 of \cite{brediesSuperpositionInhomogeneous2022}.
The normalization $h(0)=1$ is obtained by absorbing the initial weight
into the representing measure. 
\end{proof}

Proposition \ref{prop:BCF-controlled} provides a Lagrangian representation of an unbalanced
transportation. In contrast to the classical superposition principle, where an elementary trajectory contributes
the unit mass \(\delta_{\gamma(t)}\), here each characteristic contributes the weighted mass
\(h(t)\delta_{\gamma(t)}\). The spatial trajectory \(\gamma\) describes transport, whereas \(h\) records the
cumulative creation or destruction of mass along that trajectory. This is the reason the natural path space
for the present problem is the augmented space of weighted characteristics.

\subsection{Relaxation by weighted arcs}

The use of Young-measure-valued controls follows the classical relaxed-control
framework \cite{young1969lectures,warga1972optimalControl}.

Let $\mathcal{P}(U)$ denote the set of Borel probability measures
on $U$. A \emph{relaxed weighted arc} is a triple $a=(\gamma,h,\nu)$
such that 
\[
\gamma\in AC([0,T];X),\quad h\in AC([0,T];\mathbb{R}_{>0}),\quad\nu_{t}\in\mathcal{P}(U)
\]
with $t\mapsto\nu_{t}$ weakly measurable, $h(0)=1$, and 
\begin{align}
\dot{\gamma}(t) & =\int_{U}F(t,\gamma(t),u)d\nu_{t}(u),\label{eq:relaxed-arc-x}\\
\dot{h}(t) & =h(t)\int_{U}G(t,\gamma(t),u)d\nu_{t}(u)\label{eq:relaxed-arc-h}
\end{align}
for a.e. $t$. Denote the set of all such arcs by $\mathcal{A}$.
It is convenient to identify the relaxed control $\nu$ with the Young
measure $dt\,d\nu_{t}(u)$ on $[0,T]\times U$ and equip these measures
with the weak-{*} topology. Since $U$ is compact, the relaxed control
set is compact. On $\mathcal{A}$ we use the product of this topology
and uniform convergence of $(\gamma,h)$. The bounds on $F$ and $G$
over the compact set $[0,T]\times X\times U$ imply uniform Lipschitz
bounds for $\gamma$ and $h$; hence the admissible relaxed-arc set
is compact by Arzelà-Ascoli and weak-{*} compactness of Young measures.
The relaxed dynamics are closed under this convergence.

Define the arc cost 
\begin{equation}
c(a):=\int^{T}_{0}h(t)\int_{U}L(t,\gamma(t),u)d\nu_{t}(u)dt.\label{eq:relaxed-arc-cost}
\end{equation}

A feasible relaxed path measure is a finite nonnegative measure $\Pi\in\mathcal{M}_{+}(\mathcal{A})$
satisfying 
\begin{equation}
(e_{0}\circ\pi_{1})_{\#}\Pi=\mu_{0}\label{eq:relaxed-initial}
\end{equation}
and 
\begin{equation}
\int_{\mathcal{A}}h(T)\psi(\gamma(T))d\Pi(a)=\int_{X}\psi d\mu_{T},\quad\forall\psi\in C(X).\label{eq:relaxed-terminal}
\end{equation}

\medskip{}
 \noindent\textbf{Problem P.} Let $\mu_{0}\in\mathcal{M}_{+}(K),\,\mu_{T}\in\mathcal{M}_{+}(X)$,
define
\begin{equation}
V_{P}(\mu_{0},\mu_{T}):=\inf_{\Pi}\int_{\mathcal{A}}c(a)d\Pi(a),\label{eq:VP}
\end{equation}
where the infimum is over all feasible relaxed path measures.

Once the admissible arc set $\mathcal{A}$ is fixed, both the objective
and the endpoint constraints are linear in the positive measure $\Pi$.
Thus Problem P is an infinite-dimensional linear program over measures.
It is a relaxation of Problem E because ordinary controls correspond
to Dirac relaxed controls $\nu_{t}=\delta_{u(t)}$. 
\begin{thm}[Relaxation is tight]
\label{thm:relax_tight} Under Assumption A, 
\begin{equation}
V_{E}(\mu_{0},\mu_{T})=V_{P}(\mu_{0},\mu_{T}).\label{eq:relax_tight}
\end{equation}
\end{thm}
The proof of the theorem is provided in the Appendix.

\begin{remark}
Theorem \ref{thm:relax_tight} shows that the path-space formulation is an exact convexification rather than 
merely a lower
relaxation of the Eulerian problem E. The key structural assumptions are the affinity of \(F\) and \(G\) in the
control input, the convexity of \(U\), and the convexity of \(L\). Without these convexity properties,
Problem P should in general be interpreted as a relaxed-control problem and a relaxation gap may
occur.
\end{remark}

\section{Duality }\label{sec:duality}

Since Problem P is a linear program, we can, in principle derive its
Kantorovich dual. Define the Hamiltonian 
\begin{equation}
\begin{aligned}H(t,x,p,p_{0}) & =\sup_{u\in U}\bigl\{ p^{\top}F(t,x,u)\\
 & \qquad\qquad+p_{0}G(t,x,u)-L(t,x,u)\bigr\}.
\end{aligned}
\label{eq:Hamiltonian}
\end{equation}
The scalar variable $p_{0}$ reflects the reaction term in \eqref{sys:PDE}.

\subsection{The smooth HJB dual problem}

For $v\in C^{1}([0,T]\times X)$, define 
\begin{equation}
J_{D}(v):=\int_{X}v(T,x)\,d\mu_{T}(x)-\int_{X}v(0,x)\,d\mu_{0}(x).\label{eq:JD}
\end{equation}

\medskip{}
 \noindent\textbf{Problem D.} Let $\mu_{0}\in\mathcal{M}_{+}(K),\,\mu_{T}\in\mathcal{M}_{+}(X)$,
define
\begin{equation}
V_{D}(\mu_{0},\mu_{T}):=\sup_{v}J_{D}(v),\label{eq:VD}
\end{equation}
where the supremum is over $v\in C^{1}([0,T]\times X)$ satisfying
\begin{equation}
\partial_{t}v+H(t,x,\nabla v,v)\le0\label{eq:HJB-subsolution}
\end{equation}
for $(t,x)\in[0,T]\times X$. Here $C^{1}([0,T]\times X)$ is understood
as the restriction to $[0,T]\times X$ of a $C^{1}$ function defined
on a neighborhood. 

It is easy to see that
\begin{equation}
V_{D}(\mu_{0},\mu_{T})\le V_{E}(\mu_{0},\mu_{T})=V_{P}(\mu_{0},\mu_{T}).\label{eq:weak-dual}
\end{equation}
In fact, for any feasible Eulerian pair and any $v$ satisfying \eqref{eq:HJB-subsolution},
the definition of $H$ gives 
\[
\partial_{t}v+\nabla v^{\top}F(t,x,u)+G(t,x,u)v\le L(t,x,u).
\]
Integrating this inequality against the continuity equation gives
$J_{D}(v)\le J_{E}(\mu,u)$. Take the supremum in $v$ and then the
infimum in $(\mu,u)$, and use Theorem \ref{thm:relax_tight}. 

In the next subsection, we proceed to prove the strong duality between
Probems P and D. 

\subsection{Endpoint dual of the relaxed problem}

Define the endpoint value functional 
\begin{equation}
\mathcal{C}_{\mu_{0}}(\nu):=V_{P}(\mu_{0},\nu),\quad\nu\in\mathcal{M}_{+}(X),\label{eq:endpoint-functional}
\end{equation}
with value $+\infty$ when the endpoint is not reachable with finite
cost. 
\begin{prop}[Convexity and lower semicontinuity]
\label{prop:endpoint-lsc} The endpoint functional $\mathcal{C}_{\mu_{0}}$
is proper whenever at least one finite-cost endpoint is reachable,
convex, and weak-{*} lower semicontinuous on $\mathcal{M}(X)$. 
\end{prop}
\begin{proof}
Convexity follows immediately from the linearity of Problem P in the
path measure $\Pi$. For lower semicontinuity, let $\nu_{k}\stackrel{*}{\rightharpoonup}\nu$
and choose nearly optimal path measures $\Pi_{k}$ for $\mathcal{C}_{\mu_{0}}(\nu_{k})$.
Their total mass is fixed by the initial constraint and the relaxed-arc
space is compact in the topology introduced above. Hence, after extraction,
$\Pi_{k}\stackrel{*}{\rightharpoonup}\Pi$. The initial and weighted
terminal constraints pass to the limit, while continuity of the data
and lower semicontinuity of the integral cost give 
\[
\int_{\mathcal{A}}c(a)d\Pi(a)\le\liminf_{k\to\infty}\int_{\mathcal{A}}c(a)d\Pi_{k}(a).
\]
Thus $\Pi$ is feasible for the endpoint $\nu$ and the desired lower
semicontinuity follows. 
\end{proof}

For $\psi\in C(X)$, define the terminal-reward value function 
\begin{equation}
W_{\psi}(t,x)=\sup_{u(\cdot)}h(T)\psi(\gamma(T))-\int^{T}_{t}h(s)L(s,\gamma,u)ds\label{eq:terminal-reward-value}
\end{equation}
subject to
\[
\dot{\gamma}=F(s,\gamma,u),\quad\dot{h}=G(s,\gamma,u)h,\quad\gamma(t)=x,\quad h(t)=1.
\]
The supremum may equivalently be taken over relaxed controls, because
$F$ and $G$ are affine in $u$ and $L$ is convex in $u$. 
\begin{prop}[Conjugate representation]
\label{prop:conjugate} For every $\psi\in C(X)$, 
\begin{equation}
\mathcal{C}^{*}_{\mu_{0}}(\psi)=\int_{X}W_{\psi}(0,\cdot)d\mu_{0}.\label{eq:conjugate-value}
\end{equation}
Consequently, 
\begin{equation}
V_{P}(\mu_{0},\mu_{T})=\sup_{\psi\in C(X)}\left\{ \int_{X}\psi d\mu_{T}-\int_{X}W_{\psi}(0,\cdot)d\mu_{0}\right\} .\label{eq:nonsmooth-dual}
\end{equation}
\end{prop}

The representation has a simple control-theoretic interpretation. Taking the convex conjugate removes the
prescribed terminal marginal and replaces it by the terminal reward \(\psi\). Once the terminal measure is
free, the path-space optimization separates with respect to the initial particles, and the resulting
one-particle problem is precisely \(W_\psi\). Thus \(W_\psi(0,\cdot)\) plays the role of a dynamic nonlinear
analogue of the $c$-transform \cite{villaniOptimalTransport2009}.

Equation \eqref{eq:nonsmooth-dual} is already the natural dual of
the path-space linear program. The remaining difficulty is to identify
this dual value with Problem D, whose admissible functions are $C^{1}$
HJB subsolutions.

\subsection{From the value function to smooth subsolutions}

The next lemma separates the regularity of the value function from
the smooth-subsolution approximation used in the duality theorem. 
\begin{lem}[Regularity of the terminal-reward value]
\label{lem:W-regularity} If $\psi$ is Lipschitz on $X$, then $W_{\psi}$
is Lipschitz on every compact cylinder containing the relevant controlled
trajectories. The dynamic programming principle identifies $W_{\psi}$
as the viscosity solution of 
\begin{equation}
\partial_{t}W_{\psi}+H(t,x,\nabla W_{\psi},W_{\psi})=0,\quad W_{\psi}(T,\cdot)=\psi.\label{eq:W-HJB}
\end{equation}
In particular, the HJB equality holds at almost every differentiability
point of $W_{\psi}$. 
\end{lem}
\begin{proof}
Use the same control for trajectories issued from two initial states.
The local Lipschitz bounds for $F$, $G$, $L$, together with Gronwall's
inequality and the boundedness of the trajectories and weights, yield
a spatial Lipschitz estimate uniform over controls. The dynamic programming
principle then gives the corresponding Lipschitz estimate in time.
The viscosity characterization is the standard dynamic programming
argument. 
\end{proof}

\begin{remark}
The Lipschitz regularity does not imply that $W_{\psi}$
is a classical $C^{1}$ solution. Indeed, nonsmoothness may develop
through crossing optimal characteristics. This is precisely why the
passage to the smooth dual class requires a separate approximation
theorem; compare \cite{bernardOptimalMassTransportation2007} and
the classical convex duality constructions of \cite{flemingConvexDualityApproach1989,vinterConvexDualityNonlinear1993}.
\end{remark}

\begin{thm}[Smooth $C^{1}$-subsolution approximation]
\label{thm:smooth-subsolution} Let $\psi$ be the restriction to
$X$ of a $C^{1}$ function defined on a neighborhood of $X$. For
every $\varepsilon>0$, there exists $v_{\varepsilon}\in C^{1}([0,T]\times X)$
such that 
\begin{equation}
\partial_{t}v_{\varepsilon}+H(t,x,\nabla v_{\varepsilon},v_{\varepsilon})\le0,\label{eq:smooth-HJB}
\end{equation}
\begin{equation}
v_{\varepsilon}(T,x)\ge\psi(x),\quad x\in X,\label{eq:smooth-terminal}
\end{equation}
and for all $x\in{\rm supp}\mu_{0}$,
\begin{equation}
W_{\psi}(0,x)\le v_{\varepsilon}(0,x)\le W_{\psi}(0,x)+\varepsilon,.\label{eq:smooth-initial}
\end{equation}
\end{thm}

\begin{remark}
The approximation in Theorem \ref{thm:smooth-subsolution} is deliberately one-sided. Strong duality does not
require a global smooth approximation of the viscosity solution \(W_\psi\). It is enough to construct a smooth
HJB subsolution that majorizes the prescribed terminal reward and approximates \(W_\psi\) from above at the
initial states carrying mass. These are exactly the two properties entering the dual objective. The
construction therefore avoids imposing classical solvability of the HJB equation.
\end{remark}

We can now identify the smooth HJB dual. 
\begin{thm}[Strong duality]
\label{thm:strong-duality} Under Assumption A, 
\begin{equation}
V_{E}(\mu_{0},\mu_{T})=V_{P}(\mu_{0},\mu_{T})=V_{D}(\mu_{0},\mu_{T}).\label{eq:main-duality}
\end{equation}
\end{thm}
\begin{proof}
The equality $V_{E}=V_{P}$ follows from
Theorem~\ref{thm:relax_tight}, and weak
duality \eqref{eq:weak-dual} gives $V_{D}\le V_{P}$. In
\eqref{eq:nonsmooth-dual}, the supremum may be restricted to $C^{1}$
terminal data by density and the continuity of the terminal-reward
value with respect to the terminal datum. Fix such a $\psi$ and choose
$v_{\varepsilon}$ from Theorem \ref{thm:smooth-subsolution}. Then
$v_{\varepsilon}$ is feasible for Problem D and 
\begin{align*}
J_{D}(v_{\varepsilon}) & \ge\int_{X}\psi d\mu_{T}-\int_{X}W_{\psi}(0,x)d\mu_{0}(x)-\varepsilon|\mu_{0}|.
\end{align*}
Since $\epsilon$ is arbitrary, we get $V_{D}\ge V_{P}$. 
\end{proof}

Theorem \ref{thm:strong-duality} identifies three equivalent descriptions of the same problem, the original
Eulerian control
problem, a convex linear program over weighted trajectories, and a smooth HJB dual problem. In particular, the
passage to relaxed path measures introduces zero relaxation gap, while restricting the dual problem to
\(C^1\) subsolutions introduces no duality gap.

\section{Reachability}

\label{sec:reachability}

So far, we have not studied the feasibility of the dynamical UOT problem.
Reachability for continuity equations and measure-valued control systems
has been studied in several balanced settings \cite{cavagnariGeneralizedControlSystems2018,cavagnariAveragedTimeOptimal2018,cavagnariMarigondaAttainability2021,caponigroTransportMaps2025}.
With the presence of the reaction term, the reachability problem shall
be treated separately

Thanks to the duality theory, we are now able to characterize the
measures that can be reached from a prescribed initial measure at
time $T$. Throughout this section, $\mu_{0}\in\mathcal{M}_{+}(K)$
is fixed, where $X$ is the compact localization set introduced in
Section~\ref{sec:problem}. 

\subsection{The reachable set}
\begin{defn}
\label{def:reachable-set} Let $\mu_{0}\in\mathcal{M}_{+}(K)$, we
say that $\mu_{T}\in\mathcal{M}_{+}(X)$ is reachable from $\mu_{0}$
if there exists an admissible Eulerian pair $(\mu,u)$ satsifying
\eqref{sys:PDE} and $\mu(0)=\mu_{0}$, $\mu(T)=\mu_{T}$. Denote
the set of all reachable measures as $\mathcal{R}_{T}(\mu_{0})$.
\end{defn}

By Theorem~\ref{thm:relax_tight}, the same set is obtained from
the relaxed weighted-arc formulation. More precisely, let $\mathcal{A}$
denote the set of relaxed weighted arcs introduced in Section~\ref{sec:relaxation},
and define 
\begin{equation}
\mathfrak{P}(\mu_{0}):=\left\{ \Pi\in\mathcal{M}_{+}(\mathcal{A}):(e_{0}\circ\pi_{1})_{\#}\Pi=\mu_{0}\right\} .\label{eq:path-measures-fixed-initial}
\end{equation}
For $\Pi\in\mathfrak{P}(\mu_{0})$, define its terminal weighted marginal
$\mathcal{E}_{T}(\Pi)\in\mathcal{M}_{+}(X)$ by 
\begin{equation}
\int_{X}\psi(x)d\mathcal{E}_{T}(\Pi)(x):=\int_{\mathcal{A}}h(T)\psi(\gamma(T))d\Pi(a),\label{eq:terminal-map}
\end{equation}
for $\psi\in C(X)$. Then 
\begin{equation}
\mathcal{R}_{T}(\mu_{0})=\left\{ \mathcal{E}_{T}(\Pi):\Pi\in\mathfrak{P}(\mu_{0})\right\} .\label{eq:reachable-set-path}
\end{equation}

We first note some elementary geometric properties of the reachable
set. 
\begin{prop}
\label{prop:reachable-compact-convex} The set $\mathcal{R}_{T}(\mu_{0})$
is convex and narrowly compact. Moreover, if 
\begin{align}
\underline{G} & :=\min_{(t,x,u)\in[0,T]\times X\times U}G(t,x,u),\label{eq:G-bounds}\\
\overline{G} & :=\max_{(t,x,u)\in[0,T]\times X\times U}G(t,x,u),
\end{align}
then every $\mu_{T}\in\mathcal{R}_{T}(\mu_{0})$ satisfies 
\begin{equation}
e^{\underline{G}T}|\mu_{0}|\le|\mu_{T}|\le e^{\overline{G}T}|\mu_{0}|.\label{eq:reachable-mass-bounds}
\end{equation}
\end{prop}
\begin{proof}
The set $\mathfrak{P}(\mu_{0})$ is convex and the terminal map $\mathcal{E}_{T}$
defined by \eqref{eq:terminal-map} is linear. Hence \eqref{eq:reachable-set-path}
immediately implies that $\mathcal{R}_{T}(\mu_{0})$ is convex.

For a relaxed weighted arc $a=(\gamma,h,\nu)$, one has 
\[
\dot{h}(t)=h(t)\int_{U}G(t,\gamma(t),u)d\nu_{t}(u),\quad h(0)=1.
\]
Therefore, 
\[
e^{\underline{G}t}\le h(t)\le e^{\overline{G}t},\quad t\in[0,T].
\]
Using $(e_{0}\circ\pi_{1})_{\#}\Pi=\mu_{0}$, hence $|\Pi|=|\mu_{0}|$,
we obtain 
\[
|\mathcal{E}_{T}(\Pi)|=\int_{\mathcal{A}}h(T)d\Pi(a),
\]
and \eqref{eq:reachable-mass-bounds} follows.

Finally, by the compactness result for relaxed arcs established in
Section~\ref{sec:relaxation}, $\mathcal{A}$ is compact in the relaxed-arc
topology. Thus the set of positive measures on $\mathcal{A}$ with
fixed initial marginal $\mu_{0}$ is narrowly compact. Since 
$ a\mapsto h(T)\psi(\gamma(T))$
is continuous for every $\psi\in C(X)$, the map $\mathcal{E}_{T}$
is narrowly continuous. Consequently, \eqref{eq:reachable-set-path}
implies that $\mathcal{R}_{T}(\mu_{0})$ is narrowly compact. 
\end{proof}

\begin{remark}
Notice that convexity here is a property of the reachable set in the space of measures, and does not require
convexity of the pointwise reachable sets of the underlying state dynamics. The mass estimate
\eqref{eq:reachable-mass-bounds} gives, in
addition, an elementary necessary condition for reachability before any HJB test is performed.
\end{remark}

\subsection{Support function of the reachable set}

For $\psi\in C(X)$, let 
\begin{equation}
\sigma_{\mathcal{R}_{T}(\mu_{0})}(\psi):=\sup_{\mu_{T}\in\mathcal{R}_{T}(\mu_{0})}\int_{X}\psi d\mu_{T}\label{eq:reachable-support-function}
\end{equation}
denote the support function of the reachable set.

To characterize \eqref{eq:reachable-support-function}, consider the
zero-running-cost terminal-reward problem by setting $L=0$ in \eqref{eq:terminal-reward-value}.
Denote this cost as $W^{0}_{\psi}(t,x)$. Then by Proposition \ref{prop:conjugate},
\begin{equation}
\sigma_{\mathcal{R}_{T}(\mu_{0})}(\psi)=\int_{X}W^{0}_{\psi}(0,x)d\mu_{0}(x).\label{eq:reachable-support-value}
\end{equation}

Consider the zero-cost Hamiltonian 
\begin{equation}
H_{0}(t,x,p,p_{0}):=\sup_{u\in U}\left\{ p^{\top}F(t,x,u)+p_{0}G(t,x,u)\right\} ,\label{eq:reachability-Hamiltonian}
\end{equation}
and the class of smooth HJB subsolutions 
\begin{equation}
\mathcal{S}_{0}:=\left\{ \varphi\in C^{1}([0,T]\times X):\partial_{t}\varphi+H_{0}(t,x,\nabla\varphi,\varphi)\le0\right\} .\label{eq:reachability-subsolutions}
\end{equation}

The smooth-subsolution approximation developed in Theorem~\ref{thm:smooth-subsolution}
applies directly with $L\equiv0$. 
\begin{cor}
\label{cor:reachable-support-HJB} For every $\psi\in C(X)$, 
\begin{equation}
\sigma_{\mathcal{R}_{T}(\mu_{0})}(\psi)=\inf_{\substack{\varphi\in\mathcal{S}_{0}\\
\varphi(T,\cdot)\ge\psi
}
}\int_{X}\varphi(0,x)d\mu_{0}(x).\label{eq:reachable-support-HJB}
\end{equation}
\end{cor}

Hence the support function of the entire reachable set can be computed without explicitly constructing
reachable measures. It is determined by smooth HJB supersolutions of the terminal payoff. Since a closed
convex set is uniquely determined by its support function, this representation leads directly to the
reachability criterion below.

\subsection{A dual characterization of reachability}

We can now characterize the reachable set entirely in terms of smooth
HJB subsolutions. 
\begin{thm}
\label{thm:reachability-characterization} Let $\mu_{T}\in\mathcal{M}_{+}(X)$.
Then $\mu_{T}\in\mathcal{R}_{T}(\mu_{0})$ if and only if 
\begin{equation}
\int_{X}\varphi(T,x)d\mu_{T}(x)\le\int_{X}\varphi(0,x)d\mu_{0}(x),\quad\forall\varphi\in\mathcal{S}_{0}.\label{eq:reachability-HJB-characterization}
\end{equation}
\end{thm}
\begin{remark}
Theorem~\ref{thm:reachability-characterization} is an HJB separation
theorem for the reachable set. Every $\varphi\in\mathcal S_{0}$ gives
a necessary inequality satisfied by all reachable terminal measures.
Hence a single violation of \eqref{eq:reachability-HJB-characterization}
provides a smooth certificate of non-reachability. Conversely, because
the reachable set is convex and narrowly compact and its support
function is given by Corollary~\ref{cor:reachable-support-HJB}, these
certificates are complete: every unreachable target is separated by at
least one smooth HJB subsolution. Thus reachability is the zero-running-cost
counterpart of the duality result developed in the preceding section.
\end{remark}

\begin{example}
\label{ex:reachability-certificate}
Consider the one-dimensional system on $X=[-1,1]$ with $T=1$ and
$U=[-1,1]$,
\[
F(x,u)=u,\qquad G(x,u)=u,
\]
and initial measure $\mu_{0}=\delta_{0}$. Transport and reaction are
therefore generated by exactly the same scalar actuator. The zero-cost
Hamiltonian is
\[
H_{0}(p,p_{0})=\sup_{|u|\le1}u(p+p_{0})=|p+p_{0}|.
\]
The two smooth functions $\varphi_{\pm}(t,x)=\pm e^{-x}$ satisfy
\[
\partial_{t}\varphi_{\pm}
+H_{0}(\partial_{x}\varphi_{\pm},\varphi_{\pm})=0,
\]
so both belong to $\mathcal S_{0}$. Consequently, every reachable
terminal measure must satisfy both inequalities in
\eqref{eq:reachability-HJB-characterization}, and hence
\begin{equation}
\int_{X}e^{-x}d\mu_{T}(x)
=\int_{X}e^{-x}d\mu_{0}(x)=1.
\label{eq:certificate-invariant}
\end{equation}
For the candidate target $\mu_{T}=\delta_{1}$, however,
\[
\int_{X}\varphi_{-}(1,x)d\mu_{T}(x)
=-e^{-1}>-1
=\int_{X}\varphi_{-}(0,x)d\mu_{0}(x),
\]
which violates \eqref{eq:reachability-HJB-characterization} and
certifies that $\delta_{1}$ is unreachable. 
\end{example}

\subsection{Sufficient conditions from the augmented dynamics }

\label{subsec:reachability-augmented}

We first provide a sufficient condition for reachability in terms
of the finite-dimensional characteristic dynamics. Let $\mu_{0},\mu_{1}\in\mathcal{M}_{+}(X)$
be nonzero and denote 
\[
M_{0}:=|\mu_{0}|,\quad M_{1}:=|\mu_{1}|.
\]
Recall that a weighted characteristic satisfies 
\[
\dot{x}(t)=F(t,x(t),u(t)),\quad\dot{h}(t)=G(t,x(t),u(t))h(t)
\]
and $h(0)=1$. Since $h(t)>0$, introducing $z(t):=\log h(t)$ gives
the augmented control system 
\begin{equation}
\begin{aligned}\dot{x}(t) & =F(t,x(t),u(t)),\\
\dot{z}(t) & =G(t,x(t),u(t)).
\end{aligned}
\label{eq:augmented-system}
\end{equation}

For $(x,z_{0}),(y,z_{1})\in X\times\mathbb{R}$, we write 
\[
(x,z_{0})\leadsto_{T}(y,z_{1})
\]
if there exists an admissible control for \eqref{eq:augmented-system}
steering $(x,z_{0})$ to $(y,z_{1})$ in time $T$. 
\begin{prop}
\label{prop:augmented-sufficient} Let $\mu_{0}\in\mathcal{M}_{+}(K),\,\mu_{1}\in\mathcal{M}_{+}(X)$
be nonzero. Suppose that 
\begin{equation}
(x,0)\leadsto_{T}\left(y,\log\frac{M_{1}}{M_{0}}\right)\label{eq:uniform-augmented-reachability}
\end{equation}
for every $x\in{\rm supp}\mu_{0},\;y\in{\rm supp}\mu_{1}$. Then $\mu_{T}\in\mathcal{R}_{T}(\mu_{0})$. 
\end{prop}
\begin{cor}
\label{cor:augmented-global} Suppose that the augmented system \eqref{eq:augmented-system}
is fixed-time controllable in the following sense: for every $x,y\in X$
and every $\alpha$ in a prescribed interval $I\subset\mathbb{R}$,
\[
(x,0)\leadsto_{T}(y,\alpha).
\]
Then any two nonzero measures $\mu_{0},\mu_{1}\in\mathcal{M}_{+}(X)$
satisfying $\log\frac{|\mu_{1}|}{|\mu_{0}|}\in I$ are mutually connected
by an admissible measure evolution from $\mu_{0}$ to $\mu_{1}$ in
time $T$. 
\end{cor}
\begin{remark} Proposition \ref{prop:augmented-sufficient} is only a
sufficient condition. Measure reachability may hold even if \eqref{eq:uniform-augmented-reachability}
fails, since the relaxed path-space formulation permits the mass issuing
from a point to be split among several admissible characteristics.
The theorem is useful because it reduces feasibility of the measure-valued
problem to a standard controllability property of the finite-dimensional
augmented system \eqref{eq:augmented-system}. \end{remark}

\subsection{A special case: balanced transport}

\label{subsec:balanced-reachability}

We now specialize to the balanced case $G(t,x,u)\equiv0$. The transport
equation becomes 
\begin{equation}
\partial_{t}\mu+\nabla\cdot\bigl(F(t,x,u)\mu\bigr)=0,\label{eq:balanced-continuity}
\end{equation}
and every admissible evolution preserves total mass.

The zero-cost Hamiltonian reduces to 
\begin{equation}
\widetilde{H}_{0}(t,x,p):=\sup_{u\in U}p^{\top}F(t,x,u).\label{eq:balanced-zero-Hamiltonian}
\end{equation}
Define the class 
\begin{equation}
\widetilde{\mathcal{S}}_{0}:=\left\{ \varphi\in C^{1}:\partial_{t}\varphi+\widetilde{H}_{0}(t,x,\partial_{x}\varphi)\le0\right\} .\label{eq:balanced-subsolutions}
\end{equation}

The following result is an immediate zero-running-cost consequence
of the duality theorem. 
\begin{prop}
\label{prop:balanced-reachability} Let $\mu_{0}\in\mathcal{M}_{+}(K),\,\mu_{T}\in\mathcal{M}_{+}(X)$.
Under the assumptions of Theorem~\ref{thm:strong-duality}, $\mu_{T}\in\mathcal{R}_{T}(\mu_{0})$
if and only if for every $\varphi\in\widetilde{\mathcal{S}}_{0}$,
\begin{equation}
\int_{X}\varphi(T,x)d\mu_{1}(x)\le\int_{X}\varphi(0,x)d\mu_{0}(x).\label{eq:balanced-HJB-reachability}
\end{equation}
\end{prop}
\begin{proof}
Consider Problem E with running cost identically zero, $L\equiv0$.
Its value with prescribed endpoints $(\mu_{0},\mu_{T})$ is 
\[
V_{0}(\mu_{0},\mu_{1})=\begin{cases}
0, & \mu_{T}\in\mathcal{R}_{T}(\mu_{0}),\\[1mm]
+\infty, & \mu_{T}\notin\mathcal{R}_{T}(\mu_{0}).
\end{cases}
\]
By Theorem~\ref{thm:strong-duality}, 
\begin{equation}
V_{0}(\mu_{0},\mu_{T})=\sup_{\varphi\in\widetilde{\mathcal{S}}_{0}}\left\{ \int_{X}\varphi(T,\cdot)d\mu_{1}-\int_{X}\varphi(0,\cdot)d\mu_{0}\right\} .\label{eq:balanced-zero-cost-dual}
\end{equation}

If $\mu_{T}$ is reachable, then $V_{0}(\mu_{0},\mu_{T})=0$. Since
$\varphi\equiv0$ is feasible in \eqref{eq:balanced-zero-cost-dual},
the supremum is at least zero. Therefore every feasible $\varphi$
must satisfy \eqref{eq:balanced-HJB-reachability}.

Conversely, suppose \eqref{eq:balanced-HJB-reachability} holds for
every $\varphi\in\widetilde{\mathcal{S}}_{0}$. Then the right-hand
side of \eqref{eq:balanced-zero-cost-dual} is at most zero. Since
$\varphi\equiv0$ is feasible, it is exactly zero. Hence $V_{0}(\mu_{0},\mu_{T})=0$,
and therefore $\mu_{T}\in\mathcal{R}_{T}(\mu_{0})$. 
\end{proof}

The preceding dual characterization can be expressed directly in terms
of the reachable sets of the underlying finite-dimensional system
\begin{equation}
\dot{x}(t)=F(t,x(t),u(t)).\label{eq:underlying-ODE}
\end{equation}
For $x\in X$, define its fixed-time reachable set by $\mathcal{R}^{{\rm ode}}_{T}(x)$.
For $\psi\in C(X)$, consider the terminal-reward value function $V^{0}_{\psi}(t,x):=\sup_{u(\cdot)}\psi(x(T))$.
At the initial time, 
\begin{equation}
V^{0}_{\psi}(0,x)=\sup_{y\in\mathcal{R}^{{\rm ode}}_{T}(x)}\psi(y).\label{eq:value-ODE-reachable}
\end{equation}

\begin{cor}
\label{cor:balanced-ODE-reachability} Let $\mu_{0}\in\mathcal{M}_{+}(K),\,\mu_{T}\in\mathcal{M}_{+}(X)$.
Then $\mu_{T}\in\mathcal{R}_{T}(\mu_{0})$ if and only if 
\begin{equation}
\int_{X}\psi(y)d\mu_{1}(y)\le\int_{X}\sup_{y\in\mathcal{R}^{{\rm ode}}_{T}(x)}\psi(y)d\mu_{0}(x),\quad\forall\psi\in C(X).\label{eq:balanced-ODE-dual}
\end{equation}
\end{cor}
\begin{proof}
For $\psi\in C(X)$, let 
\[
\sigma_{\mathcal{R}_{T}(\mu_{0})}(\psi):=\sup_{\nu\in\mathcal{R}_{T}(\mu_{0})}\int_{X}\psi d\nu
\]
be the support function of the measure reachable set. The free-terminal
version of the zero-cost path-space problem gives 
\begin{equation}\label{eq:balanced-support-ODE}
	\begin{aligned}
		\sigma_{\mathcal{R}_{T}(\mu_{0})}(\psi) 
		& =\int_{X}V^{0}_{\psi}(0,x)d\mu_{0}(x). \\
		& =\int_{X}\sup_{y\in\mathcal{R}^{{\rm ode}}_{T}(x)}\psi(y)d\mu_{0}(x).
	\end{aligned}
\end{equation}
If $\mu_{T}\in\mathcal{R}_{T}(\mu_{0})$, then by definition of the
support function, 
\[
\int_{X}\psi d\mu_{T}\le\sigma_{\mathcal{R}_{T}(\mu_{0})}(\psi),
\]
which, together with \eqref{eq:balanced-support-ODE}, yields \eqref{eq:balanced-ODE-dual}.

Conversely, suppose \eqref{eq:balanced-ODE-dual} holds for every
$\psi\in C(X)$. Then 
\[
\int_{X}\psi d\mu_{T}\le\sigma_{\mathcal{R}_{T}(\mu_{0})}(\psi),\quad\forall\psi\in C(X).
\]
Since the relaxed path-space representation implies that $\mathcal{R}_{T}(\mu_{0})$
is a convex narrowly closed subset of $\mathcal{M}(X)$, the Hahn-Banach
separation theorem implies $\mu_{T}\in\mathcal{R}_{T}(\mu_{0})$.
\end{proof}
\begin{cor}
\label{cor:balanced-global-controllability} 
If the system \eqref{eq:underlying-ODE} is controllable for every $x\in {\rm supp}\mu_0$, i.e., 
\begin{equation}
\mathcal{R}^{{\rm ode}}_{T}(x)=X,\quad\forall x\in{\rm supp}\mu_{0},\label{eq:global-fixed-time-controllability}
\end{equation}
then 
\begin{equation}
\mathcal{R}_{T}(\mu_{0})=\left\{ \nu\in\mathcal{M}_{+}(X):|\nu|=|\mu_{0}|\right\} .\label{eq:all-balanced-measures-reachable}
\end{equation}
\end{cor}
\begin{proof}
Under \eqref{eq:global-fixed-time-controllability}, 
\[
\sup_{y\in\mathcal{R}^{{\rm ode}}_{T}(x)}\psi(y)=\max_{y\in X}\psi(y)
\]
for every $x\in X$. Thus, for any positive measure $\mu_{1}$ satisfying
$|\mu_{1}|=|\mu_{0}|$, 
\[
\int_{X}\psi\,d\mu_{1}\le|\mu_{1}|\max_{X}\psi=\int_{X}\sup_{y\in\mathcal{R}^{{\rm ode}}_{T}(x)}\psi(y)\,d\mu_{0}(x).
\]
Corollary \ref{cor:balanced-ODE-reachability} therefore implies $\mu_{T}\in\mathcal{R}_{T}(\mu_{0}).$ 
\end{proof}
\begin{remark} 
	Corollary \ref{cor:balanced-ODE-reachability} and Corollary
	\ref{cor:balanced-global-controllability} provide
direct bridges between classical finite-dimensional controllability
and measure-valued reachability. The pointwise sets $\mathcal{R}^{{\rm ode}}_{T}(x)$
completely determine the support function of the reachable set of
measures through \eqref{eq:balanced-support-ODE}. In particular,
global fixed-time controllability of the underlying ODE implies reachability
of every positive target measure having the same total mass. This
result has already been obtained in \cite{elamvazhuthiDynamicalOptimalTransport2023}.
But here we have followed a different approach, namely, through duality
theory. \end{remark}

\section{Numerics }\label{sec:numerics}
The dual formulation of Section~\ref{sec:duality} provides a convenient
starting point for numerical computation. Recall that 
\begin{equation}
V_{D}(\mu_{0},\mu_{T})=\sup_{v}\left\{ J_{D}(v):\partial_{t}v+H(t,x,\nabla v,v)\le0\right\} ,\label{eq:numerical-dual}
\end{equation}
where 
\[
J_{D}(v)=\int_{X}v(T,x)d\mu_{T}(x)-\int_{X}v(0,x)d\mu_{0}(x).
\]
Since $H(t,x,\cdot,\cdot)$ is convex, the HJB inequality defines
a pointwise convex constraint. This suggests an augmented-Lagrangian
splitting analogous to those used in dynamic optimal transport \cite{papadakisOptimalTransportProximal2014,nishinoBenamouBrenier2025}.

\subsection{Dual augmented-Lagrangian scheme}\label{subsec:dual-alg}

Define 
\begin{align*}
K_{H}(t,x) & :=\left\{ (a,p,p_{0}):a+H(t,x,p,p_{0})\le0\right\} ,\\
Av & :=(\partial_{t}v,\nabla v,v).
\end{align*}
Then the dual problem is equivalently written as 
\begin{equation}
\min_{v,q}\left\{ -J_{D}(v)+I_{K_{H}}(q)\right\} \text{ subject to }q=Av.\label{eq:dual-splitting-continuous}
\end{equation}

Let $A_{h}$ denote a consistent space-time discretization of $A$,
and let $J_{D,h}$ be the corresponding discretization of $J_{D}$.
Denote the optimal value of the resulting discrete dual problem by
$V_{D,h}$. For a penalty parameter $\beta>0$, the scaled ADMM iterations are
\begin{align}
v^{k+1} & =\arg\min_{v}-J_{D,h}(v)+\frac{\beta}{2}\|A_{h}v-q^{k}+y^{k}\|^{2}_{h},\label{eq:BB-v-update}\\
q^{k+1} & =\Pi_{K_{H,h}}\bigl(A_{h}v^{k+1}+y^{k}\bigr),\label{eq:BB-q-update}\\
y^{k+1} & =y^{k}+A_{h}v^{k+1}-q^{k+1}.\label{eq:BB-y-update}
\end{align}

The $v$-subproblem is quadratic. If $b_{h}$ represents the endpoint
functional $J_{D,h}$, then 
\begin{equation}
\beta A^{*}_{h}A_{h}v^{k+1}=b_{h}+\beta A^{*}_{h}(q^{k}-y^{k}).\label{eq:v-linear-system}
\end{equation}
For the standard discretization 
\[
A^{*}_{h}A_{h}=D^{*}_{t}D_{t}+\nabla^{*}_{h}\nabla_{h}+I,
\]
so each $v$-update reduces to solving a sparse linear system with
a coefficient matrix that is fixed throughout the iterations.

The nonlinear step \eqref{eq:BB-q-update} is completely pointwise,
since $K_{H,h}$ is the Cartesian product of the sets $K_{H}(t_{j},x_{\ell})$.
At each grid point one computes 
\begin{equation}
\begin{aligned}\Pi_{K_{H}}(\bar{a},\bar{p},\bar{p}_{0}) & =\arg\min_{(a,p,p_{0})\in K_{H}}[(a-\bar{a})^{2}\\
 & \qquad\qquad\qquad+|p-\bar{p}|^{2}+(p_{0}-\bar{p}_{0})^{2}].
\end{aligned}
\label{eq:projection-K-H}
\end{equation}
Hence all nonlinear projection steps can be performed independently
over the space-time grid. For the control-affine dynamics \eqref{eq:F-G},
the Hamiltonian takes the form 
\begin{equation}
H(t,x,p,p_{0})=p^{\top}f_{0}+p_{0}g_{0}+\sup_{u\in U}\left\{ c^{\top}u-L(t,x,u)\right\} ,\label{eq:H-control-affine-numerical}
\end{equation}
where $c_{i}=p^{\top}f_{i}+p_{0}g_{i}$. Thus evaluation of $H$ requires
only a finite-dimensional convex optimization over the control variable.
In particular, for 
\[
L(t,x,u)=q(t,x)+\frac{1}{2}u^{\top}R(t,x)u,\quad R(t,x)\succ0,
\]
the maximizing control is 
\begin{equation}
u^{\sharp}=\Pi^{R}_{U}(R^{-1}c),\label{eq:u-hamiltonian-projection}
\end{equation}
where $\Pi^{R}_{U}$ denotes the projection onto $U$ in the $R$-weighted
norm,
\[
\Pi^{R}_{U}(z):=\arg\min_{u\in U}\frac{1}{2}(u-z)^{\top}R(u-z).
\]
Consequently, the nonlinear part of each ADMM iteration is local and
naturally parallelizable; the only globally coupled operation is the
sparse linear system \eqref{eq:v-linear-system}.

\subsection{Recovery of primal variables and feedback}
\label{subsec:primal-recovery}

Although the algorithm is formulated entirely in terms of the HJB dual,
its multiplier has a direct primal interpretation. Let
$\lambda=(\rho,M,Q)$ and define
\begin{equation}
\sigma_{K_H}(\rho,M,Q)
:=\sup_{(a,p,p_0)\in K_H}\{\rho a+M^{\top}p+Qp_0\}.
\label{eq:support-KH-active}
\end{equation}
Since $K_H=\{(a,p,p_0):a+H(t,x,p,p_0)\le0\}$,
\begin{equation}
\sigma_{K_H}(\rho,M,Q)=
\begin{cases}
\rho H^*\!\left(t,x,\dfrac{M}{\rho},\dfrac{Q}{\rho}\right),&\rho>0,\\[1mm]
0,&(\rho,M,Q)=0,\\
+\infty,&\text{otherwise},
\end{cases}
\label{eq:support-KH-perspective}
\end{equation}
where the conjugate is taken with respect to $(p,p_0)$. Under
Assumption A,
\begin{equation}
\label{eq:Hstar-control-active}
	\begin{aligned}
	H^*(t,x,\xi,\zeta)
	=\min\{&L(t,x,u)  :u\in U, \\
			&  F(t,x,u)=\xi, G(t,x,u)=\zeta\},
	\end{aligned}
\end{equation}
with value $+\infty$ if the constraint set is empty. Thus finiteness of
\eqref{eq:support-KH-perspective} enforces $\rho\ge0$ and, on
$\{\rho>0\}$, requires $(M/\rho,Q/\rho)$ to be generated by the same
control input.

Writing $J_{D,h}(v)=\langle b_h,v\rangle_h$, finite-dimensional Fenchel
duality gives
\begin{equation}
V_{D,h}=\min_{\lambda_h}\sigma_{K_{H,h}}(\lambda_h)
\quad\text{subject to}\quad A_h^*\lambda_h=b_h,
\label{eq:discrete-flux-primal}
\end{equation}
where $\sigma_{K_{H,h}}$ is the grid-weighted sum of
\eqref{eq:support-KH-active}. With
$\lambda_h=(\rho_h,M_h,Q_h)$, the adjoint constraint is the discrete
counterpart of
\begin{equation}
\partial_t\rho+\nabla\!\cdot M=Q,\qquad
\rho(0,\cdot)=\mu_0,\quad \rho(T,\cdot)=\mu_T.
\label{eq:flux-continuity}
\end{equation}
Hence, in the scaled ADMM iteration,
\begin{equation}
\lambda_h^k=\beta y^k=(\rho_h^k,M_h^k,Q_h^k)
\label{eq:primal-from-multiplier}
\end{equation}
recovers the discrete mass and transport-reaction fluxes at convergence.
The feedback is obtained directly from the HJB potential as
\begin{equation}
\begin{aligned}
u_h^\star(t,x)\in\arg\max_{u\in U}\{&\nabla_hv_h^\star(t,x)^\top F(t,x,u) \\ 
 +&v_h^\star(t,x)G(t,x,u) -L(t,x,u)\}.
\end{aligned}
\label{eq:feedback-from-dual}
\end{equation}
For $L=q+\frac12u^\top Ru$, this reduces to
\begin{equation}
u_h^\star=\Pi_U^R\!\left(R^{-1}c_h^\star\right),\qquad
(c_h^\star)_i=\nabla_hv_h^{\star\top}f_i+v_h^\star g_i.
\label{eq:feedback-quadratic}
\end{equation}
Thus solving the HJB dual also provides the associated discrete primal
evolution and feedback law.

\subsection{Numerical experiments: Cancer treatment} \label{subsec:numerical-experiments}

\subsubsection{Phenotype-structured cell population under treatment}
\label{subsubsec:phenotype-example}

Phenoty\-pe-structured models of cancer-cell populations have been used
to describe treatment-induced selection, drug resistance, and reversible
phenotypic adaptation \cite{chisholmDrugTolerance2015,lorzPhenotypicHeterogeneity2015,almeidaCytotoxicTherapy2019}.
In particular, Chisholm et al.~\cite{chisholmDrugTolerance2015}
introduce a two-trait description in terms of survival potential and
proliferation potential, and they model drug-induced adaptation through
phenotypic dynamics at the individual level and an advection mechanism
at the population level. Motivated by these models, we consider the
following deterministic control-affine dynamics, which may be viewed
as a simplified ODE representation of the drug-induced phenotypic
evolution,
\begin{align}
\dot{r} & =-a_{r}r+b_{r}u(1-r),\label{eq:cell-r-dynamics}\\
\dot{p} & =a_{p}(1-\chi r-p)-b_{p}u(1-r)p.\label{eq:cell-p-dynamics}\\
\dot{h} & =\bigl(\gamma p-d_{0}-k_{{\rm kill}}u(1-r)\bigr)h
\end{align}

Here $r\in[0,1]$ represents a normalized survival or drug-resistance
phenotype, $p\in[0,1]$ represents proliferative activity, and $h>0$
denotes the size of the corresponding cell subpopulation. The interpretation
of $(r,p)$ is therefore closely related to the survival-proliferation
phenotype space introduced in \cite{chisholmDrugTolerance2015}, although
the particular control-affine drift used here is a simplified deterministic
model chosen to isolate the transport-reaction mechanism of interest.

In the absence of treatment, \eqref{eq:cell-r-dynamics} models relaxation
of drug resistance toward the sensitive phenotype, whereas \eqref{eq:cell-p-dynamics}
drives proliferative activity toward the resistance-dependent level
$1-\chi r$. The factor $\chi r$ represents a fitness cost associated
with resistance. Treatment increases resistance through the term $b_{r}u(1-r)$
and suppresses proliferation through $b_{p}u(1-r)p$. Such treatment-induced
movement toward weakly proliferative and drug-tolerant phenotypes
is consistent with the phenotypic adaptation mechanisms considered
in \cite{chisholmDrugTolerance2015}. The third equation accounts
for population growth and death: $\gamma p$ represents phenotype-dependent
proliferation, $d_{0}$ basal death, and $k_{{\rm kill}}u(1-r)$ drug-induced
killing, which is stronger for drug-sensitive cells.

The system therefore has the form
\begin{equation}
\dot{x}=F(x,u),\quad\dot{h}=G(x,u)h,\label{eq:cell-characteristic-system}
\end{equation}
where 
\begin{equation}
F(r,p,u)=\begin{pmatrix}-a_{r}r\\[1mm]
a_{p}(1-\chi r-p)
\end{pmatrix}+u\begin{pmatrix}b_{r}(1-r)\\[1mm]
-b_{p}(1-r)p
\end{pmatrix},\label{eq:cell-F}
\end{equation}
and 
\begin{equation}
G(r,p,u)=\gamma p-d_{0}-k_{{\rm kill}}u(1-r).\label{eq:cell-G}
\end{equation}
The treatment intensity is denoted by 
\[
u\in U:=[0,u_{\max}].
\]
The important feature is that the same treatment $u$ simultaneously
changes the phenotype $(r,p)$ and the population weight $h$. Hence
the transport and reaction mechanisms cannot be chosen independently,
which is exactly the weighted characteristic dynamics considered in
this work.

\subsubsection{Population-level optimal transport problem.}

We now consider a population of such cells. Let $\mu_{t}\in\mathcal{M}_{+}(X)$
denote the distribution of cancer cells over phenotype space at time
$t$. The characteristic model \eqref{eq:cell-characteristic-system}
induces the controlled transport-reaction equation 
\begin{equation}
\partial_{t}\mu+\nabla_{(r,p)}\cdot\bigl(F(r,p,u)\mu\bigr)=G(r,p,u)\mu.\label{eq:cell-uot}
\end{equation}

We prescribe an initial phenotype distribution $\mu_{0}$ and a desired
terminal distribution $\mu_{T}$, with different total masses allowed.
The objective is to determine the minimum-energy treatment that steers
$\mu_{0}$ to $\mu_{T}$: 
\begin{equation}
\begin{aligned}\inf_{\mu,u}\quad & \int^{T}_{0}\int_{X}\frac{\eta}{2}u(t,r,p)^{2}d\mu_{t}(r,p)dt,\\
\mathrm{s.t.}\quad & 0\le u(t,r,p)\le u_{\max},
\end{aligned}
\label{eq:cell-UOT-problem}
\end{equation}
This is an unbalanced optimal transport problem because the total
mass is not conserved. 

\subsubsection{Simulation results}

For the numerical experiment we use $a_{r}=0.35$, $b_{r}=1.5$, $a_{p}=0.55$,
$\chi=0.6$, $b_{p}=2$, $\gamma=0.85$, $d_{0}=0.12$, $k_{{\rm kill}}=2.6$
, $u_{\max}=1.5$, $\eta=0.35$ with terminal time $T=1.5$. The initial
distribution is concentrated around a sensitive and highly proliferative
phenotype, whereas the prescribed terminal distribution has smaller
total mass and is concentrated around a more resistant and less proliferative
phenotype. More precisely, we take the truncated Gaussian 
\begin{align*}
\rho_{0}(r,p) &
=\frac{1}{Z_{0}}\exp\left[-\frac{1}{2}\left(\frac{r-0.18}{0.065}\right)^{2}-\frac{1}{2}\left(\frac{p-0.78}{0.080}\right)^{2}\right],
\end{align*}
where $Z_{0}$ is the normalization constant on $X$. The terminal
density $\rho_{T}$ is specified as
\[
	\rho_T(r,p) = \frac{0.55}{Z_T}
	\exp\left[ -\frac12 \left( \frac{r-0.56}{0.11} \right)^2 - 
		\frac12  \left( \frac{p-0.47}{0.091} \right)^2
	\right]
\]where $Z_T$ is the normalization constant on $X$.
Thus the desired transformation requires simultaneously a reduction
of approximately $45\%$ in total cell population and a substantial
change in its phenotypic composition.

\begin{figure}[h]
\centering
\includegraphics[scale=0.6]{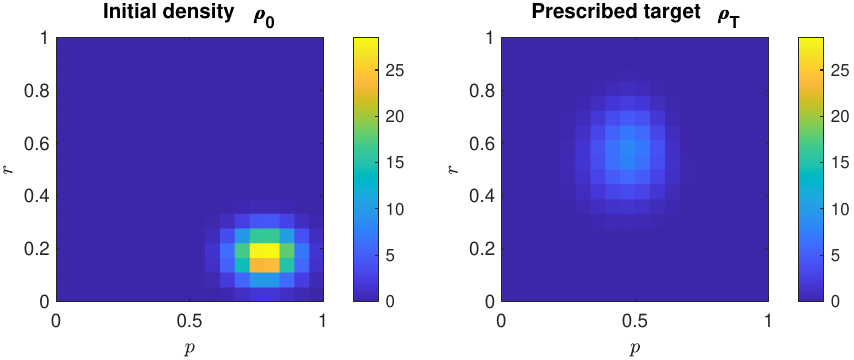}

\caption{The initial distribution and the desired distribution at the terminal
time. }\label{fig:cell-endpoints}

\end{figure}

We solve the hard-endpoint problem \eqref{eq:cell-UOT-problem} using the algorithm introduced in Sections~\ref{subsec:dual-alg}
and~\ref{subsec:primal-recovery}

\begin{figure}[h]
\centering
\includegraphics[scale=0.65]{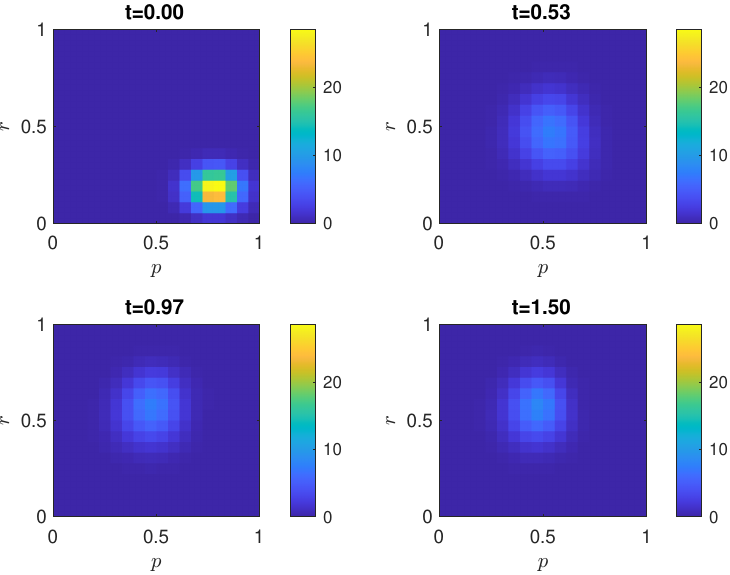}

\caption{Cell evolution. The figure shows that the control input steers the
initial distribution to the desired distribution. The $L_{2}$ error
between the controlled terminal and the desired distribution is around
$7.87\times10^{-4}$. }\label{fig:cell-evolution.}
\end{figure}

\begin{figure}[h]
\centering
\includegraphics[scale=0.33]{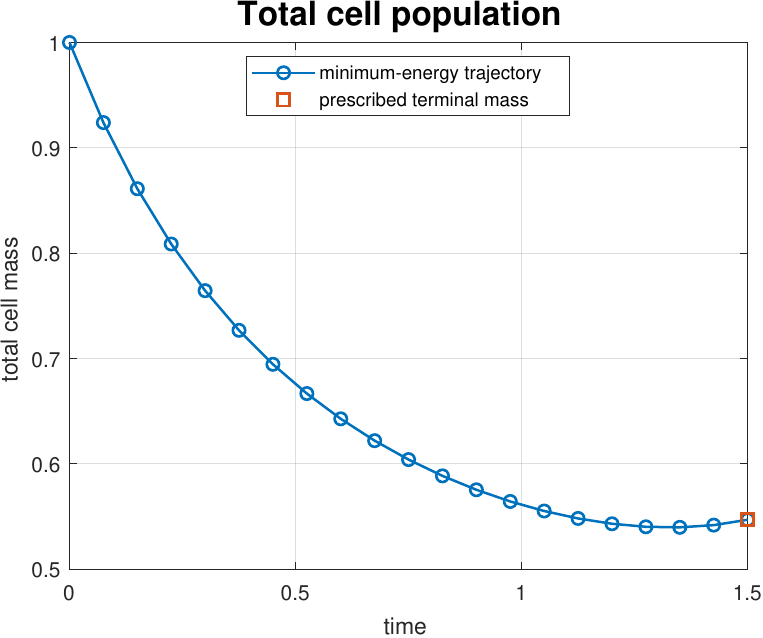}\includegraphics[scale=0.33]{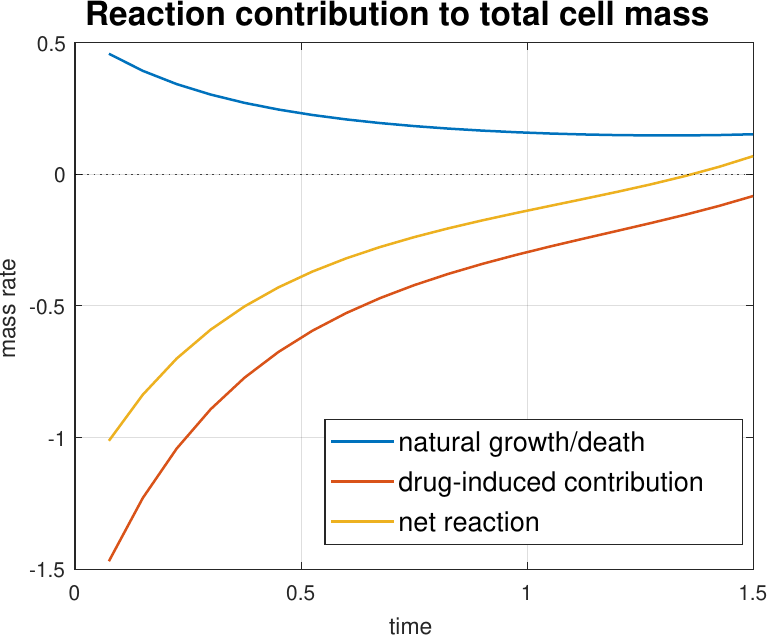}

\caption{Mass evolution of the cell population. }\label{fig:mass-evolution}
\end{figure}

Fig. \ref{fig:cell-evolution.} shows the evolution of the cell population.
From the figure, we can see that the optimal treatment transports
the population from the initially sensitive and highly proliferative
region toward a more resistant and less proliferative phenotype while
simultaneously reducing its total mass. The mean resistance increases
from $0.1805$ to $0.5606$, whereas the mean proliferative activity
decreases from $0.7794$ to $0.4670$. The total population decreases
from $M(0)=1$ to the prescribed value $M(T)=0.55$ as shown in
Fig. \ref{fig:mass-evolution}. 

\begin{figure}[h]
\centering
\includegraphics[scale=0.65]{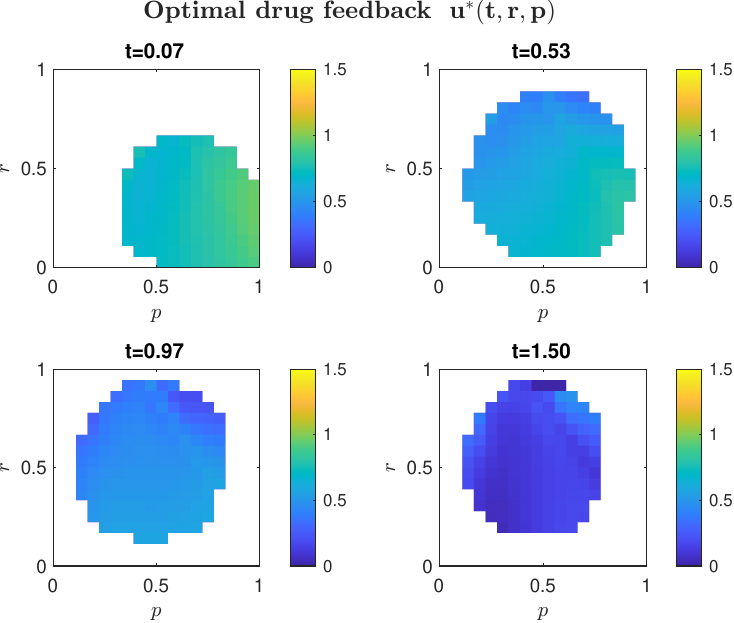}

\caption{The snapshots of the control actions. Only those values above $10^{-3}$
are shown. }\label{fig:control-actions}
\end{figure}

Fig. \ref{fig:control-actions} shows the snapshots of control action
where we display $u^{*}$ only on phenotype cells for which $\rho$
exceeds $10^{-3}$ of its instantaneous maximum; outside this occupied
region the quotient $u=w/\rho$ is physically irrelevant and may show
spurious saturated values. 

\begin{figure}[h]
\centering
\includegraphics[scale=0.33]{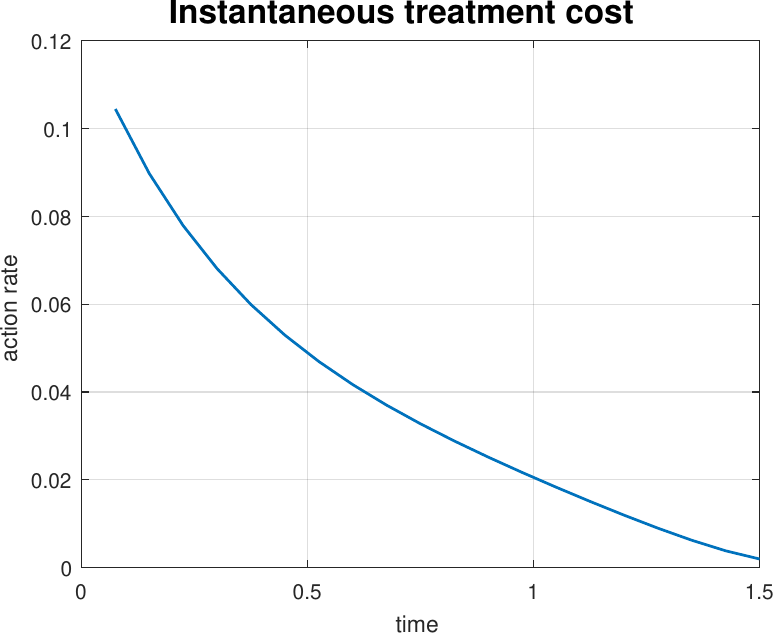}\includegraphics[scale=0.33]{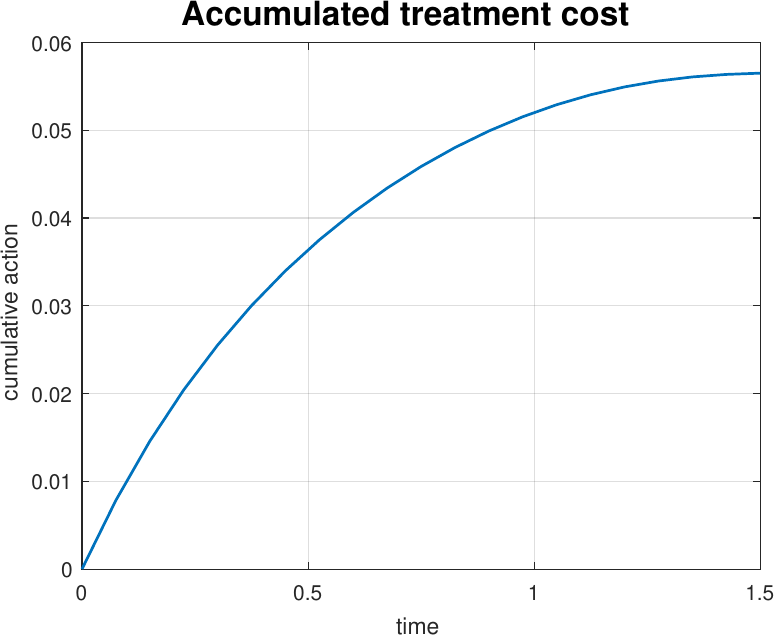}\caption{The instantaneous and accumulated treatment costs. Approximately $53\%$,
$81\%$, and $96\%$ of this action has already been spent by $t=0.375$,
$0.75$, and $1.125$, respectively. }\label{fig:action}
\end{figure}

Fig. \ref{fig:action} shows the instantaneous and accumulated treatment
cost. From the figures, we see that the treatment
is strongly front-loaded; 
indeed, approximately $53\%$, $81\%$, and $96\%$ of the action has
already been spent by $t=0.375$, $0.75$, and $1.125$, respectively.
Thus the optimizer uses relatively strong treatment in the early stage
to both kill sensitive cells and reshape the phenotype distribution,
and then progressively relaxes the treatment once the target composition
has nearly been reached. This interpretation is made explicit by the
mass balance 
\begin{equation}
\dot{M}(t)=\underbrace{\int_{X}(\gamma p-d_{0})\rho drdp}_{\mathcal{R}_{{\rm nat}}(t)}-\underbrace{\int_{X}k_{{\rm kill}}u(1-r)\rho drdp}_{\mathcal{R}_{{\rm drug}}(t)}.\label{eq:cell-mass-decomposition}
\end{equation}
The natural contribution $\mathcal{R}_{{\rm nat}}$ remains positive
over the entire horizon, whereas the drug-induced contribution is
initially strongly negative. Their sum is negative over most of the
treatment interval, causing a rapid decrease in cell burden. Near
the end of the horizon, the treatment has been sufficiently reduced
that the net reaction changes sign. Accordingly, the total population
reaches a minimum of approximately $0.5396$ at $t\approx1.35$ and
then undergoes a small amount of regrowth in order to meet the prescribed
terminal mass. This behavior is not a numerical artifact: it reflects
the fact that \eqref{eq:cell-UOT-problem} seeks minimum energy subject
to an exact terminal distribution, rather than maximal cell killing.

\section{Conclusion}

This paper developed a control-theoretic framework for dynamical unbalanced optimal transport in which
transport and mass variation are coupled through a common control input. Relaxation and duality form the two
central components of the framework. The former lifts the original nonlinear problem to a linear one 
over a larger path space, while the latter brings the problem back to the state space through a
Hamilton-Jacobi-Bellman characterization. Thus, relaxation unveils the hidden convex structure of
the problem, whereas duality provides a tractable representation that is directly connected to
control-theoretic objects.

The dual formulation further serves as the basis for two subsequent developments. On the theoretical side, it
leads to a characterization of reachability in terms of HJB inequalities and separating
certificates. On the computational side, it gives rise to numerical algorithms posed entirely on the dual
variables.

The present analysis relies on control-affine dynamics, a convex control set and running cost. Extending the
framework to noncompact state spaces, stochastic or diffusive dynamics, and scalable grid-free methods for
higher-dimensional problems provides several directions for further study.

\section*{Appendix}
\begin{proof}[Proof of Proposition \ref{prop:conjugate}]
Extend $\mathcal{C}_{\mu_{0}}$ to $\mathcal{M}(X)$ by setting $\mathcal{C}_{\mu_{0}}(\nu)=+\infty$
whenever $\nu$ is not a nonnegative reachable endpoint with finite
cost. By definition of the convex conjugate, 
\begin{equation}
\mathcal{C}^{*}_{\mu_{0}}(\psi)=\sup_{\nu\in\mathcal{M}(X)}\left\{ \int_{X}\psi d\nu-\mathcal{C}_{\mu_{0}}(\nu)\right\} .\label{eq:conjugate-definition}
\end{equation}
For a fixed reachable endpoint $\nu$, the definition of Problem P
gives 
\[
\mathcal{C}_{\mu_{0}}(\nu)=\inf_{\Pi\in\Gamma(\mu_{0},\nu)}\int_{\mathcal{A}}c(a)d\Pi(a),
\]
where $\Gamma(\mu_{0},\nu)$ denotes the set of relaxed path measures
satisfying \eqref{eq:relaxed-initial} and \eqref{eq:relaxed-terminal}
with terminal measure $\nu$. Hence 
\begin{align}
 & \int_{X}\psi\,d\nu-\mathcal{C}_{\mu_{0}}(\nu)\nonumber \\
 & \qquad=\sup_{\Pi\in\Gamma(\mu_{0},\nu)}\left\{ \int_{X}\psi d\nu-\int_{\mathcal{A}}c(a)d\Pi(a)\right\} .\label{eq:conjugate-fixed-endpoint}
\end{align}
Using the terminal constraint \eqref{eq:relaxed-terminal}, 
\[
\int_{X}\psi d\nu=\int_{\mathcal{A}}h(T)\psi(\gamma(T))d\Pi(a).
\]
Therefore, after taking the supremum over the free terminal measure
$\nu$, \eqref{eq:conjugate-definition} becomes 
\begin{equation}
\mathcal{C}^{*}_{\mu_{0}}(\psi)=\sup_{\substack{\Pi\in\mathcal{M}_{+}(\mathcal{A})\\
(\gamma(0))_{\#}\Pi=\mu_{0}
}
}\int_{\mathcal{A}}\mathcal{R}_{\psi}(a)\,d\Pi(a),\label{eq:conjugate-path}
\end{equation}
where 
\begin{equation}
\mathcal{R}_{\psi}(a):=h(T)\psi(\gamma(T))-c(a).\label{eq:path-reward}
\end{equation}
Thus the terminal constraint has disappeared: the conjugate is a free-terminal
relaxed optimal-control problem with prescribed initial measure.

For $x\in\operatorname{supp}\mu_{0}$, let 
\[
\mathcal{A}(x):=\{a=(\gamma,h,\nu)\in\mathcal{A}:\gamma(0)=x\}.
\]
Since ordinary controls are included among relaxed controls through
$\nu_{t}=\delta_{u(t)}$, while every relaxed control can be replaced
by its barycenter without changing $(\gamma,h)$ and without increasing
the running cost, affinity of $F,G$ and convexity of $L$ imply 
\begin{equation}
W_{\psi}(0,x)=\sup_{a\in\mathcal{A}(x)}\mathcal{R}_{\psi}(a).\label{eq:W-relaxed-arcs}
\end{equation}

Let $\Pi$ be arbitrary in \eqref{eq:conjugate-path}. For $\Pi$-a.e.
arc $a=(\gamma,h,\nu)$, 
\[
\mathcal{R}_{\psi}(a)\le W_{\psi}(0,\gamma(0))
\]
by \eqref{eq:W-relaxed-arcs}. Hence 
\begin{align*}
\int_{\mathcal{A}}\mathcal{R}_{\psi}(a)\,d\Pi(a) & \le\int_{\mathcal{A}}W_{\psi}(0,\gamma(0))\,d\Pi(a)\\
 & =\int_{X}W_{\psi}(0,x)\,d\mu_{0}(x),
\end{align*}
where the last equality follows from $(e_{0}\circ\pi_{1})_{\#}\Pi=\mu_{0}$.
Taking the supremum over $\Pi$ in \eqref{eq:conjugate-path} gives
\begin{equation}
\mathcal{C}^{*}_{\mu_{0}}(\psi)\le\int_{X}W_{\psi}(0,x)\,d\mu_{0}(x).\label{eq:conjugate-upper}
\end{equation}

For the reverse inequality, fix $\varepsilon>0$. The relaxed-arc
space $\mathcal{A}$ is compact, the map $a\mapsto\gamma(0)$ is continuous,
and $\mathcal{R}_{\psi}$ is continuous under the topology introduced
above. Consequently, the multifunction 
\(
x\mapsto\mathcal{A}(x)
\)
has compact values and closed graph. Standard measurable-selection
results \cite[Theorem 18.19]{AliprantisBorder2006} therefore yield a Borel map 
\[
x\longmapsto a^{\varepsilon}_{x}=(\gamma^{\varepsilon}_{x},h^{\varepsilon}_{x},\nu^{\varepsilon}_{x})\in\mathcal{A}(x)
\]
such that 
\begin{equation}
\mathcal{R}_{\psi}(a^{\varepsilon}_{x})\ge W_{\psi}(0,x)-\varepsilon,\quad x\in\operatorname{supp}\mu_{0}.\label{eq:epsilon-optimal-selection}
\end{equation}
Define 
\begin{equation}
\Pi_{\varepsilon}:=\bigl(x\mapsto a^{\varepsilon}_{x}\bigr)_{\#}\mu_{0}.\label{eq:Pi-epsilon-conjugate}
\end{equation}
Since $\gamma^{\varepsilon}_{x}(0)=x$, 
\(
(e_{0}\circ\pi_{1})_{\#}\Pi_{\varepsilon}=\mu_{0},
\)
so $\Pi_{\varepsilon}$ is admissible in \eqref{eq:conjugate-path}.
Therefore, 
\begin{align*}
\mathcal{C}^{*}_{\mu_{0}}(\psi) & \ge\int_{\mathcal{A}}\mathcal{R}_{\psi}(a)d\Pi_{\varepsilon}(a)\\
 & =\int_{X}\mathcal{R}_{\psi}(a^{\varepsilon}_{x})d\mu_{0}(x)\\
 & \ge\int_{X}W_{\psi}(0,x)d\mu_{0}(x)-\varepsilon|\mu_{0}|.
\end{align*}
Letting $\varepsilon\downarrow0$ yields 
\begin{equation}
\mathcal{C}^{*}_{\mu_{0}}(\psi)\ge\int_{X}W_{\psi}(0,x)d\mu_{0}(x).\label{eq:conjugate-lower}
\end{equation}
Together with \eqref{eq:conjugate-upper}, this proves \eqref{eq:conjugate-value}.

Finally, Proposition \ref{prop:endpoint-lsc} shows that $\mathcal{C}_{\mu_{0}}$
is proper, convex, and weak-{*} lower semicontinuous. Hence the Fenchel-Moreau
theorem on the dual pair $(\mathcal{M}(X),C(X))$ gives 
\[
\mathcal{C}_{\mu_{0}}(\mu_{T})=\sup_{\psi\in C(X)}\left\{ \int_{X}\psi d\mu_{T}-\mathcal{C}^{*}_{\mu_{0}}(\psi)\right\} .
\]
Substituting \eqref{eq:conjugate-value} yields \eqref{eq:nonsmooth-dual}. 
\end{proof}

\bigskip
\begin{proof}[Proof of Theorem \ref{thm:relax_tight}]
Let $(\mu,u)$ be feasible for Problem E. By Proposition \ref{prop:BCF-controlled},
there exists a measure $\sigma$ on $\mathcal{A}_{0}$ satisfying
\eqref{eq:path-initial} - \eqref{eq:path-h}. For $\sigma$-a.e.
$z=(\gamma,h)$, define the relaxed control $\xi_{t}\in\mathcal{P}(U)$
as 
\begin{equation}
\xi^{z}_{t}:=\delta_{u(t,\gamma(t))},\quad\text{for a.e. }t\in[0,T],\label{eq:Dirac-relaxed-control}
\end{equation}
and the corresponding relaxed weighted arc 
\[
\Theta(z):=(\gamma,h,t\mapsto\xi^{z}_{t})\in\mathcal{A}.
\]
Define 
\begin{equation}
\Pi:=(\Theta)_{\#}\sigma.\label{eq:Pi-from-sigma}
\end{equation}
Since $\xi^{z}_{t}$ is a Dirac measure, \eqref{eq:path-x} and \eqref{eq:path-h}
become 
\begin{align*}
\dot{\gamma}(t) & =\int_{U}F(t,\gamma(t),v)d\xi^{z}_{t}(v),\\
\dot{h}(t) & =h(t)\int_{U}G(t,\gamma(t),v)d\xi^{z}_{t}(v),
\end{align*}
so $\Theta(z)\in\mathcal{A}$ for $\sigma$-a.e. $z$. Moreover, 
\[
(e_{0}\circ\pi_{1})_{\#}\Pi=(e_{0}\circ\pi_{1})_{\#}\sigma=\mu_{0},
\]
and, for every $\psi\in C(X)$, 
\begin{align*}
\int_{\mathcal{A}}h(T)\psi(\gamma(T))d\Pi(a) & =\int h(T)\psi(\gamma(T))d\sigma(\gamma,h)\\
 & =\int_{X}\psi(x)d\mu_{T}(x),
\end{align*}
where the last equality follows from \eqref{eq:weighted-superposition}.
Hence $\Pi$ is feasible for Problem P.

Finally, 
\begin{align*}
\int_{\mathcal{A}}c(a)d\Pi(a) & =\int\int^{T}_{0}h(t)L\bigl(t,\gamma(t),u(t,\gamma(t))\bigr)dtd\sigma(\gamma,h)\\
 & =\int^{T}_{0}\int_{X}L(t,x,u(t,x))d\mu_{t}(x)dt\\
 & =J_{E}(\mu,u),
\end{align*}
where the second equality follows from \eqref{eq:weighted-superposition}
and Fubini's theorem. Therefore $V_{P}\le V_{E}$.

Conversely, let $\Pi$ be feasible for Problem P and define 
\begin{equation}
\mu_{t}:=\int_{\mathcal{A}}h(t)\delta_{\gamma(t)}d\Pi(a).\label{eq:mu-from-Pi}
\end{equation}
Introduce a measure $\eta$ as follows 
\begin{align}
\int\zeta(t,x,u)d\eta=\int_{\mathcal{A}}\int^{T}_{0}h(t)\int_{U}\zeta(t,\gamma(t),u)d\nu_{t}(u)dtd\Pi(a)\label{eq:occupation-relaxed}
\end{align}
for bounded Borel $\zeta$. Its $(t,x)$-marginal is $dt\mu_{t}(dx)$,
so disintegration \cite[Theorem 5.3.1]{ambrosioGigliSavare2008} gives
\[
d\eta(t,x,u)=d\lambda_{t,x}(u)d\mu_{t}(x)dt
\]
for a Borel family $\lambda_{t,x}\in\mathcal{P}(U)$. Define 
\begin{equation}
\bar{u}(t,x):=\int_{U}ud\lambda_{t,x}(u).\label{eq:barycentric-control}
\end{equation}
Since $U$ is convex, $\bar{u}(t,x)\in U$. Integrating the chain
rule along the arcs and using \eqref{eq:relaxed-initial} - \eqref{eq:relaxed-terminal}
shows that $\mu$ satisfies \eqref{sys:PDE} with the averaged coefficients
$\int_{U}Fd\lambda_{t,x}$ and $\int_{U}Gd\lambda_{t,x}$. By affinity
in $u$, 
\begin{align*}
 & \int_{U}F(t,x,u)d\lambda_{t,x}(u)=F(t,x,\bar{u}(t,x)),\\
 & \int_{U}G(t,x,u)d\lambda_{t,x}(u)=G(t,x,\bar{u}(t,x)).
\end{align*}
The representation \eqref{eq:mu-from-Pi} and dominated convergence
show that $t\mapsto\mu_{t}$ is narrowly continuous. Moreover, all
arcs remain in $X$, so the boundedness of $F$ and $G$ on $[0,T]\times X\times U$
implies \eqref{eq:HK-energy}. Thus $(\mu,\bar{u})$ is feasible for
Problem E. Finally, Jensen's inequality and the convexity of $L$
give 
\begin{align*}
J_{E}(\mu,\bar{u}) & \le\int^{T}_{0}\int_{X}\int_{U}L(t,x,u)d\lambda_{t,x}(u)d\mu_{t}(x)dt\\
 & =\int_{\mathcal{A}}c(a)d\Pi(a).
\end{align*}
Taking infima yields $V_{E}\le V_{P}$. 
\end{proof}

\bigskip
\begin{proof}[Proof of Theorem \ref{thm:smooth-subsolution}]
Let $\widetilde{\psi}\in C^{1}(\mathbb{R}^{n})$ be an extension
of $\psi$ which agrees with $\psi$ on a neighborhood of $X$. Since
$X$ is compact, such an extension can be chosen with compact support
after multiplication by a suitable smooth cutoff.

We first enlarge slightly the space-time region on which the value
function is considered. Fix $\eta>0$ and define 
\[
X_{\eta}:=\{x\in\mathbb{R}^{n}:\operatorname{dist}(x,X)\le\eta\}.
\]
Extend $F$, $G$, and $L$ to $[-\eta,T]\times\mathbb{R}^{n}\times U$
by freezing them at $t=0$ for negative times, namely, 
$\widetilde{N}(t,x,u) = N(t_+,x,u)$, where $t_{+}:=\max\{t,0\}$ and $N$ can be $F,G$ and $L$.
Denote the corresponding Hamiltonian
by $\widetilde{H}$. Thus 
\begin{align*}
\widetilde{H}(t,x,p,p_{0}) & =\sup_{u\in U}\{p^{\top}\widetilde{F}(t,x,u)\\
 & \qquad\qquad+p_{0}\widetilde{G}(t,x,u)-\widetilde{L}(t,x,u)\}.
\end{align*}
Clearly, 
\( \widetilde{H}(t,x,p,p_{0})=H(t,x,p,p_{0}),\) for $t\in[0,T]$.

Consider the terminal-reward value function associated with the extended
data and terminal reward $\widetilde{\psi}$, denoted by $W_{\widetilde{\psi}}$,
on $[-\eta,T]\times X_{\eta}$. By the same trajectory estimates as
in Lemma \ref{lem:W-regularity}, $W_{\widetilde{\psi}}$ is Lipschitz
on this compact cylinder. Moreover, by the dynamic programming principle,
\begin{equation}
\partial_{t}W_{\widetilde{\psi}}+\widetilde{H}\bigl(t,x,\partial_{x}W_{\widetilde{\psi}},W_{\widetilde{\psi}}\bigr)=0\label{eq:extended-W-HJB}
\end{equation}
at almost every differentiability point.

Notice that the extension does not alter the value at the initial
states relevant to the transport problem. Indeed, if $x\in\operatorname{supp}\mu_{0}$,
every admissible trajectory issued from $x$ remains in $X$ by the
localization established above. Therefore, 
\begin{equation}
W_{\psi}(0,x)=W_{\widetilde{\psi}}(0,x),\quad x\in\operatorname{supp}\mu_{0}.\label{eq:extension-independence}
\end{equation}
Also, 
\[
W_{\widetilde{\psi}}(T,x)=\widetilde{\psi}(x)=\psi(x),\qquad x\in X.
\]

We now regularize $W_{\widetilde{\psi}}$. Let 
\[
\rho\in C^{\infty}_{c}((0,1)\times B_{1}(0))
\]
be non-negative and satisfy 
\(
\int\rho(\tau,y)d\tau dy=1.
\)
For $\delta>0$, set 
\[
	\rho_{\delta}(\tau,y)=\delta^{-(n+1)}\rho\left(\frac{\tau}{\delta},\frac{y}{\delta}\right).
\]
For $\delta<\eta$, define 
\begin{equation}
W^{\delta}(t,x):=\int\rho_{\delta}(\tau,y)W_{\widetilde{\psi}}(t-\tau,x-y)d\tau dy,\label{eq:one-sided-mollification}
\end{equation}
for $(t,x)\in[0,T]\times X$, The support of $\rho_{\delta}$ is contained
in $(0,\delta)\times B_{\delta}(0)$, so only values with $t-\tau\le t\le T$
occur in \eqref{eq:one-sided-mollification}. The extension to negative
times is therefore sufficient to define the convolution also near
$t=0$.

Since $W_{\widetilde{\psi}}$ is Lipschitz, $W^{\delta}$ is smooth
on a neighborhood of $[0,T]\times X$, and 
\begin{equation}
W^{\delta}\to W_{\widetilde{\psi}}\quad\text{uniformly on }[0,T]\times X\label{eq:mollification-uniform}
\end{equation}
as $\delta\to0$.

We next estimate the HJB residual. Let $L_{W}$ be a Lipschitz constant
of $W_{\psi}$ on the enlarged cylinder and let $M_{W}:=\sup|W_{\psi}|$
there. Since $\widetilde{H}$ is continuous and the set 
\[
[-\eta,T]\times X_{\eta}\times\overline{B}_{L_{W}}(0)\times[-M_{W},M_{W}]
\]
is compact, $\widetilde{H}$ is uniformly continuous on this set.
Hence there exists a modulus $\omega(\delta)\downarrow0$ such that
\begin{align}
\big|\widetilde{H}(t,x,p,p_{0})-\widetilde{H}(t-\tau,x-y,p,p_{0})\big|\le\omega(\delta)\label{eq:H-modulus}
\end{align}
whenever 
\[
0\le t\le T,\quad x\in X,\quad0\le\tau\le\delta,\quad|y|\le\delta,
\]
and $|p|\le L_{W},\:|p_{0}|\le M_{W}$. 

Because $\widetilde{H}(t,x,\cdot,\cdot)$ is convex in $(p,p_{0})$,
Jensen's inequality gives 
\begin{align*}
 & \partial_{t}W^{\delta}(t,x)+H\bigl(t,x,\nabla W^{\delta}(t,x),W^{\delta}(t,x)\bigr)\\
 & \le\int\rho_{\delta}(\tau,y)[\partial_{t}W_{\psi}(t-\tau,x-y)\\
 & \qquad+\widetilde{H}\bigl(t,x,\nabla W_{\psi}(t-\tau,x-y),W_{\psi}(t-\tau,x-y)\bigr)]d\tau dy.
\end{align*}
Using \eqref{eq:H-modulus}, this is bounded above by 
\begin{align*}
 & \int\rho_{\delta}(\tau,y)[\partial_{t}W_{\psi}(t-\tau,x-y)+\Psi(t,\tau,x,y)]d\tau dy+\omega(\delta)
\end{align*}
where 
\begin{align*}
 & \Psi(t,\tau,x,y)\\
= & \widetilde{H}\bigl(t-\tau,x-y,\nabla W_{\psi}(t-\tau,x-y),W_{\psi}(t-\tau,x-y)\bigr).
\end{align*}
The integrand vanishes almost everywhere by \eqref{eq:extended-W-HJB}.
Therefore 
\begin{equation}
\partial_{t}W^{\delta}+H(t,x,\nabla W^{\delta},W^{\delta})\le\omega(\delta)\label{eq:mollified-HJB-residual}
\end{equation}
on $[0,T]\times X$.

It remains to remove the residual and ensure the terminal inequality.
Define 
\begin{equation}
\kappa(\delta):=\sup_{x\in X}\bigl(\psi(x)-W^{\delta}(T,x)\bigr)_{+}.\label{eq:kappa-delta}
\end{equation}
By \eqref{eq:mollification-uniform} and $W_{\psi}(T,\cdot)=\psi$
on $X$, $\kappa(\delta)\to0$. Since $X\times[0,T]\times U$ is compact,
there exists 
\[
R:=\sup_{(t,x,u)\in[0,T]\times X\times U}|G(t,x,u)|<\infty.
\]
Consequently, the Hamiltonian is $R$-Lipschitz in its last argument:
\begin{equation}
|H(t,x,p,p_{0})-H(t,x,p,q_{0})|\le R|p_{0}-q_{0}|.\label{eq:H-p0-Lipschitz}
\end{equation}

Define 
\begin{equation}
v_{\delta}(t,x):=W^{\delta}(t,x)+a_{\delta}(t)\label{eq:v-delta}
\end{equation}
where $a_{\delta}(t)\ge0$ is a function to be determined. Using \eqref{eq:mollified-HJB-residual},
\eqref{eq:H-p0-Lipschitz}, and the fact that $a_{\delta}\ge0$, we
obtain 
\begin{align*}
\partial_{t}v_{\delta}+H(t,x,\partial_{x}v_{\delta},v_{\delta})\le & \partial_{t}W^{\delta}+H(t,x,\partial_{x}W^{\delta},W^{\delta})\\
 & \qquad\qquad\qquad+a_{\delta}'+Ra_{\delta}\\
\le & \omega(\delta)+a_{\delta}'+Ra_{\delta}.
\end{align*}
Clearly, we can choose $a_{\delta}$ such that $a_{\delta}+Ra_{\delta}=0$,
with boundary condition $a_{\delta}(T)=\kappa(\delta)$. By this choice,
$v_{\delta}$ satisfies \eqref{eq:smooth-HJB}.

At the terminal time, 
\[
v_{\delta}(T,x)=W^{\delta}(T,x)+\kappa(\delta)\ge\psi(x),\quad x\in X,
\]
so \eqref{eq:smooth-terminal} also holds.

It remains to establish \eqref{eq:smooth-initial}. First, we show
the lower bound. Let $x\in\operatorname{supp}\mu_{0}$ and let $u(\cdot)$
be an arbitrary admissible open-loop control, with corresponding trajectory
$\gamma(\cdot)$ and weight $h(\cdot)$ satisfying 
\[
\dot{\gamma}=F(t,\gamma,u),\quad\dot{h}=G(t,\gamma,u)h,\quad\gamma(0)=x,\quad h(0)=1.
\]
By construction of $X$, the trajectory remains in $X$. Since $v_{\delta}$
satisfies \eqref{eq:smooth-HJB}, for every $u\in U$, 
\[
\partial_{t}v_{\delta}+\nabla v^{\top}_{\delta}F(t,x,u)+G(t,x,u)v_{\delta}\le L(t,x,u).
\]
Therefore, along the characteristic, 
\begin{align*}
\frac{d}{dt}\bigl[h(t)v_{\delta}(t,\gamma(t))\bigr] & =h(t)\Bigl(\partial_{t}v_{\delta}+\nabla v^{\top}_{\delta}F+Gv_{\delta}\Bigr)\\
 & \le h(t)L(t,\gamma(t),u(t)).
\end{align*}
Integrating from $0$ to $T$ gives 
\[
v_{\delta}(0,x)\ge h(T)v_{\delta}(T,\gamma(T))-\int^{T}_{0}h(t)L(t,\gamma(t),u(t))dt.
\]
Using the terminal majorization, 
\[
v_{\delta}(0,x)\ge h(T)\psi(\gamma(T))-\int^{T}_{0}h(t)L(t,\gamma(t),u(t))dt.
\]
Taking the supremum over all controls yields 
\begin{equation}
v_{\delta}(0,x)\ge W_{\psi}(0,x).\label{eq:v-lower-W}
\end{equation}

Finally, by \eqref{eq:mollification-uniform}, 
\[
\sup_{x\in\operatorname{supp}\mu_{0}}|W^{\delta}(0,x)-W_{\psi}(0,x)|\to0,
\]
while $a_{\delta}(0)\to0$ as $\delta\to0$. Hence, for $\delta>0$
sufficiently small, 
\[
W^{\delta}(0,x)+a_{\delta}(0)\le W_{\psi}(0,x)+\varepsilon
\]
uniformly on $\operatorname{supp}\mu_{0}$. Combining this with \eqref{eq:v-lower-W}
gives 
\[
W_{\psi}(0,x)\le v_{\delta}(0,x)\le W_{\psi}(0,x)+\varepsilon,\quad x\in\operatorname{supp}\mu_{0}.
\]
Choosing such a $\delta$ and setting $v_{\varepsilon}:=v_{\delta}$
proves the theorem. 
\end{proof}

\bigskip
\begin{proof}[Proof of Corollary \ref{cor:reachable-support-HJB}]
We first assume that $\psi$ is the restriction to $X$ of a $C^{1}$
function. Applying Theorem~\ref{thm:smooth-subsolution} with $L\equiv0$
gives, for every $\varepsilon>0$, a function $\varphi_{\varepsilon}\in\mathcal{S}_{0}$
such that 
\(
\varphi_{\varepsilon}(T,\cdot)\ge\psi
\)
and 
\[
W^{0}_{\psi}(0,x)\le\varphi_{\varepsilon}(0,x)\le W^{0}_{\psi}(0,x)+\varepsilon,\quad x\in\operatorname{supp}\mu_{0}.
\]
Consequently, 
\[
\inf_{\substack{\varphi\in\mathcal{S}_{0}\\
\varphi(T)\ge\psi
}
}\int_{X}\varphi(0,x)d\mu_{0}(x)\le\int_{X}W^{0}_{\psi}(0,x)\,d\mu_{0}(x)+\varepsilon|\mu_{0}|.
\]
Letting $\varepsilon\downarrow0$ and using \eqref{eq:reachable-support-value}
gives one inequality in \eqref{eq:reachable-support-HJB}.

For the converse, let $\varphi\in\mathcal{S}_{0}$ satisfy $\varphi(T,\cdot)\ge\psi$.
Along every controlled weighted characteristic, 
\begin{align*}
\frac{d}{dt}\left[h(t)\varphi(t,\gamma(t))\right] & =h(t)\left(\partial_{t}\varphi+\nabla\varphi^{\top}F+G\varphi\right)\le0.
\end{align*}
Therefore $\varphi(0,x)\ge h(T)\psi(\gamma(T))$ for every admissible
trajectory starting from $x$. Taking the supremum over controls gives
$\varphi(0,x)\ge W^{0}_{\psi}(0,x)$. Hence 
\[
\int_{X}\varphi(0,x)\,d\mu_{0}(x)\ge\int_{X}W^{0}_{\psi}(0,x)\,d\mu_{0}(x),
\]
and \eqref{eq:reachable-support-value} yields the reverse inequality.

For a general $\psi\in C(X)$, choose $\psi_{k}\in C^{1}(X)$ such
that 
\[
\psi_{k}\ge\psi,\quad\|\psi_{k}-\psi\|_{C(X)}\to0.
\]
Since $|G|\le R$ on $[0,T]\times X\times U$, 
\begin{equation}
\sup_{x\in X}\left|W^{0}_{\psi_{1}}(0,x)-W^{0}_{\psi_{2}}(0,x)\right|\le e^{RT}\|\psi_{1}-\psi_{2}\|_{C(X)}.\label{eq:zero-cost-value-continuity}
\end{equation}
Thus 
\(
\sigma_{\mathcal{R}_{T}(\mu_{0})}(\psi_{k})\to\sigma_{\mathcal{R}_{T}(\mu_{0})}(\psi),
\)
and passage to the limit proves \eqref{eq:reachable-support-HJB}. 
\end{proof}

\bigskip
\begin{proof}[Proof of Proposition \ref{prop:balanced-reachability}]
Let $\widehat{\mu}_{0}:=\frac{\mu_{0}}{M_{0}},\;\widehat{\mu}_{1}:=\frac{\mu_{1}}{M_{1}}$,
and choose an arbitrary coupling $\pi\in\Gamma(\widehat{\mu}_{0},\widehat{\mu}_{1})$.
Set $\alpha:=\log\frac{M_{1}}{M_{0}}$. By assumption, for every $(x,y)\in{\rm supp}\pi$
there exists an admissible trajectory of the augmented system satisfying
$x(0)=x$, $x(T)=y$, $z(0)=0$, $z(T)=\alpha$. Equivalently, for
the corresponding weighted characteristic, $h(0)=1$, $h(T)=e^{\alpha}=\frac{M_{1}}{M_{0}}$.

Using the compact relaxed-arc formulation developed above and a measurable
selection of such trajectories, let 
\[
(x,y)\mapsto a_{x,y}=(\gamma_{x,y},h_{x,y},\nu_{x,y})
\]
be a Borel selection of relaxed weighted arcs satisfying 
\[
\gamma_{x,y}(0)=x,\quad\gamma_{x,y}(T)=y,\quad h_{x,y}(T)=\frac{M_{1}}{M_{0}}.
\]
Define 
\begin{equation}
\Pi:=M_{0}\bigl((x,y)\mapsto a_{x,y}\bigr)_{\#}\pi.\label{eq:reachability-Pi}
\end{equation}
Then, for every $\varphi\in C(X)$, 
\begin{align*}
\int_{\mathcal{A}}\varphi(\gamma(0))d\Pi & =M_{0}\int_{X\times X}\varphi(x)d\pi(x,y)\\
 & =M_{0}\int_{X}\varphi d\widehat{\mu}_{0}=\int_{X}\varphi d\mu_{0},
\end{align*}
so the initial marginal of $\Pi$ is $\mu_{0}$.

Likewise, for every $\psi\in C(X)$, 
\begin{align*}
\int_{\mathcal{A}}h(T)\psi(\gamma(T))d\Pi & =M_{0}\int_{X\times X}\frac{M_{1}}{M_{0}}\psi(y)d\pi(x,y)\\
 & =M_{1}\int_{X}\psi d\widehat{\mu}_{1}\\
 & =\int_{X}\psi d\mu_{1}.
\end{align*}
Hence $\Pi$ is a feasible relaxed path measure joining $\mu_{0}$
to $\mu_{1}$. By the exact relaxation theorem (Theorem~\ref{thm:relax_tight}),
$\mu_{1}$ is reachable from $\mu_{0}$.
\end{proof}

\bigskip
\begin{proof}[Proof of Theorem \ref{thm:reachability-characterization}]
Suppose first that $\mu_{T}\in\mathcal{R}_{T}(\mu_{0})$ and let
$\varphi\in\mathcal{S}_{0}$. For any admissible evolution joining
$\mu_{0}$ to $\mu_{T}$, 
\[
\partial_{t}\varphi+\nabla\varphi^{\top}F(t,x,u)+G(t,x,u)\varphi\le0.
\]
Using the weak formulation of \eqref{sys:PDE} therefore gives 
\[
\int_{X}\varphi(T,x)d\mu_{T}(x)-\int_{X}\varphi(0,x)d\mu_{0}(x)\le0,
\]
which proves necessity.

Conversely, suppose that $\mu_{T}$ satisfies \eqref{eq:reachability-HJB-characterization}.
Let $\psi\in C(X)$ be arbitrary. If $\varphi\in\mathcal{S}_{0}$
satisfies $\varphi(T,\cdot)\ge\psi$, then 
\begin{align*}
\int_{X}\psi d\mu_{T} & \le\int_{X}\varphi(T,x)d\mu_{T}(x)\le\int_{X}\varphi(0,x)d\mu_{0}(x).
\end{align*}
Taking the infimum over all such $\varphi$ and applying Theorem~\ref{cor:reachable-support-HJB}
yields 
\begin{equation}
\int_{X}\psi d\mu_{T}\le\sigma_{\mathcal{R}_{T}(\mu_{0})}(\psi)=\sup_{\nu\in\mathcal{R}_{T}(\mu_{0})}\int_{X}\psi d\nu,\quad\forall\psi\in C(X).\label{eq:reachability-support-separation}
\end{equation}

By Proposition~\ref{prop:reachable-compact-convex}, $\mathcal{R}_{T}(\mu_{0})$
is a convex narrowly closed subset of $\mathcal{M}(X)$. If $\mu_{T}\notin\mathcal{R}_{T}(\mu_{0})$,
the Hahn--Banach separation theorem would provide a function $\psi\in C(X)$
such that 
\[
\int_{X}\psi d\mu_{T}>\sup_{\nu\in\mathcal{R}_{T}(\mu_{0})}\int_{X}\psi d\nu,
\]
contradicting \eqref{eq:reachability-support-separation}. Hence $\mu_{T}\in\mathcal{R}_{T}(\mu_{0})$,
which proves sufficiency. 
\end{proof}
\bibliographystyle{siam}
\bibliography{ref}
\end{document}